\documentclass[11pt]{article}
\usepackage{anysize}
\marginsize{2.0cm}{2.0cm}{2.0 cm}{2.0 cm}
\usepackage{setspace}

\usepackage{latexsym,amsfonts,amsmath}
\usepackage{graphicx}
\usepackage{color}
\usepackage{epstopdf}
\usepackage[titletoc]{appendix}

\newtheorem{condition}{Condition}[section]{\bfseries}{\itshape}

\newtheorem{assumption}{Assumption}[section]{\bfseries}{\itshape}

\newtheorem{theorem}{Theorem}[section]{\bfseries}{\itshape}

\newtheorem{corollary}{Corollary}[section]{\bfseries}{\itshape}

\newtheorem{proposition}{Proposition}[section]{\bfseries}{\itshape}

{\bfseries}{\itshape}

\newtheorem{lemma}{Lemma}[section]{\bfseries}{\itshape}

\newtheorem{remark}{Remark}[section]{\bfseries}{\itshape}

\newtheorem{definition}{Definition}[section]{\bfseries}{\itshape}

{\bfseries}{\itshape}

\begin{document}

\title{Continuous-time Markov decision processes with exponential utility}
\author{Yi
Zhang \thanks{Department of Mathematical Sciences, University of
Liverpool, Liverpool, L69 7ZL, U.K.. E-mail: yi.zhang@liv.ac.uk.}}
\date{}
\maketitle

\par\noindent{\bf Abstract:} In this paper, we consider a continuous-time Markov decision process (CTMDP) in Borel spaces, where the certainty equivalent with respect to the exponential utility of the total undiscounted cost is to be minimized. The cost rate is nonnegative. We establish the optimality equation. Under the compactness-continuity condition, we show the existence of a deterministic stationary optimal policy. We reduce the risk-sensitive CTMDP problem to an equivalent risk-sensitive discrete-time Markov decision process, which is with the same state and action spaces as the original CTMDP. In particular, the value iteration algorithm for the CTMDP problem follows from this reduction. We do not need impose any condition on the growth of the transition and cost rate in the state, and the controlled process could be explosive.
\bigskip

\bigskip


\section{Introduction}

In this paper we consider a continuous-time Markov decision process (CTMDP) in Borel state and action spaces, where a risk averse decision maker aims at minimizing the certainty equivalent of the total undiscounted cost with respect to the exponential utility. The cost rate is nonnegative. In the literature, see e.g., \cite{Bertsekas:1978,Chung:1987,Ghosh:2014,Howard:1972,Piunovski:1985}, such a problem is traditionally also referred to as risk-sensitive, or one with exponential utility or multiplicative cost. In this paper, we use these terms interchangeably. The CTMDP with a linear utility is often called risk neutral.

Ever since the pioneering paper \cite{Howard:1972} in 1972, there have been a large number of works devoted to risk-sensitive discrete-time Markov decision processes (DTMDPs), see \cite{BauerleRieder:2014,Cavazos:2000,Chung:1987,Fainberg:1982,Haskell:2015,Jaquette:1976,Jaskiewicz:2007,Jaskiewicz:2008}, to name just a few, where the term ``risk-sensitive'' is endowed with a more general meaning in the recent works \cite{BauerleRieder:2014,Haskell:2015}. The interested reader is referred to the reference list of the aforementioned works for more relevant literature. There are some significant differences between the risk-neutral and risk-sensitive problems. For example, the criterion with exponential utility is not decomposable in the sense of \cite{Fainberg:1982}, so that the corresponding convex analytic approach is more delicate and underdeveloped, c.f. \cite{Haskell:2015}. Also, consider the model with the total discounted cost; if the state and action spaces are both finite, then there is a stationary optimal policy in the risk neutral case, but not in the risk-sensitive case, c.f. \cite{Jaquette:1976,Chung:1987}. There have also been numerous works on risk-sensitive controlled diffusions, see the reference list of \cite{Ghosh:2014}.

In comparison, it is safe to say that much less work has been done on the risk-sensitive CTMDPs. To the best of our knowledge, the CTMDP with exponential utility was first considered in the less known article \cite{Piunovski:1985}, where the authors considered the problem on a finite time horizon. Only verification theorems were given, i.e., the authors discussed the consequences after one has obtained a solution to the optimality equation, provided that it satisfies certain conditions. The question of when there exits such a solution to the optimality equation was not discussed in \cite{Piunovski:1985}.  For the same problem as in \cite{Piunovski:1985}, this question was only considered in the recent papers \cite{Ghosh:2014,Wei:2016}. In \cite{Ghosh:2014} the transition rates were assumed to be bounded; and in \cite{Wei:2016} the growth of the transition rate was assumed to be bounded by some Lyapunov function. As a consequence, the controlled process in \cite{Ghosh:2014,Wei:2016} is nonexplosive under each policy. In \cite{Ghosh:2014,Wei:2016}, the cost rate was assumed to be bounded, and the arguments are based on Dynkin's formula or the Feynman--Kac formula. The author of \cite{Wei:2016} explained why it was hard to relax this boundedness condition on the cost rate if one follows the same approach as in there, see Section 7 therein. On the other hand, one should note that unbounded transition and cost rates appear in many real-life applications; as a simplest example, consider an $M/M/\infty$ queueing system with the holding cost rate being proportional to the number of enqueued customers.

By the way,  the CTMDP in \cite{Ghosh:2014,Piunovski:1985,Wei:2016} is assumed to be in a denumerable state space. In \cite{Ghosh:2014}, the infinite horizon discounted and average problems for the risk-sensitive CTMDP were also considered.

The present paper also deals with a CTMDP with exponential utility, but is rather different from the aforementioned works \cite{Ghosh:2014,Piunovski:1985,Wei:2016} in the following aspects. (a) We consider the problem of minimizing the expectation of the exponential utility of the total undiscounted cost over the infinite time horizon. In the current literature, we are not aware of other work on the infinite horizon total undiscounted cost criterion for the risk-sensitive CTMDP. (b) Our method of attack does not involve the Dynkin's formula or the Feyman-Kac formula, but is based on the reduction of the risk-sensitive CTMDP to a risk-sensitive DTMDP. As an advantage of developing this approach, we do not need bounds on the growth of the transition and cost rates, the controlled process is allowed to be explosive, and the state space is a general Borel space. Such explosive processes are related to the ``shattering into dust'' phenomenon in physics, see \cite{Wagner:2005}. Our reduction method starts with a risk-sensitive adaption of Yushkevich's method \cite{Yushkevich:1980}, which was originally proposed for risk-neutral CTMDPs, and later also developed to study piecewise deterministic Markov decision processes, see \cite{BauerleRieder:2011,Costa:2013,Davis:1993}. Compared with the case of a risk-neutral CTMDP, now both the state and action spaces of the induced DTMDP are more complicated than the original CTMDP.  A new equivalent DTMDP model with the same state and action spaces as the the original CTMDP will be induced later after further investigations. The powerful Feinberg's reduction method for risk-neutral CTMDP, see \cite{Feinberg:2004}, is not applicable because the criterion in the risk-sensitive CTMDP is not ``decomposable''.

The contributions of this paper are as follows. For a risk-sensitive CTMDP in Borel state and action spaces with the total undiscounted cost over an infinite horizon, we establish its optimality equation. Under the compactness-continuity condition, we show the existence of a deterministic stationary optimal policy. Moreover, we show that the risk-sensitive CTMDP problem is equivalent to a risk-sensitive DTMDP problem, which, as we emphasize, is with the same state and action spaces as the original CTMDP model. As a consequence of this, we also present the value iteration algorithm. Note that we only need impose rather weak conditions; the transition rate is arbitrarily unbounded, and the cost rate is arbitrarily unbounded in the state, and the controlled process can be explosive, so that there might exist no Lyapunov function.

The rest of the paper is organize as follows. We describe the controlled process and concerned optimal control problem in Section \ref{GuoZhangSec2}. In Section \ref{SICONDTMDP}, we present some results for the risk-sensitive DTMDPs that are needed for this paper.
In Section \ref{ZyExponentialSec02}, we develop Yushkevich's method to reduce the risk-sensitive CTMDP to a risk-sensitive DTMDP with more complicated state and action spaces. In Section \ref{ZyExponentialSect03}, we present and prove the main results in this paper, which is ended with a conclusion in Section \ref{ZyExponentialSect04}.

\bigskip

\par\noindent\textbf{Notations and conventions.} In what follows, ${\cal{B}}(X)$ is
the Borel $\sigma$-algebra of the topological space $X,$ $I$ stands for the indicator function, and $\delta_{\{x\}}(\cdot)$
is the Dirac measure concentrated on the singleton $\{x\},$ assumed to be measurable. A measure is $\sigma$-additive and $[0,\infty]$-valued. Below, unless stated otherwise, the term of
measurability is always understood in the Borel sense. Throughout
this article, we adopt the conventions of
\begin{eqnarray}\label{ZyExponential56}
\frac{0}{0}:=0,~0\cdot\infty:=0,~\frac{1}{0}:=+\infty,~\infty-\infty:=\infty.
\end{eqnarray}

\section{Model description and problem statement}\label{GuoZhangSec2}
The objective of this section is to describe briefly the
controlled process similarly to
\cite{Feinberg:2004,Kitaev:1986,Kitaev:1995}, and the associated
optimal control problem of interest in this paper.
\bigskip

Let $S$ be a nonempty Borel state space, $A$ be
a nonempty Borel action space, and $q$
stand for a signed kernel $q(dy|x,a)$ on ${\cal{B}}(S)$ given
$(x,a)\in S\times A$ such that
\begin{eqnarray*}
\tilde{q}(\Gamma_S|x,a):=q(\Gamma_S\setminus\{x\}|x,a)\ge 0
\end{eqnarray*}
for all $\Gamma_S\in{\cal{B}}(S).$ Throughout this article we assume
that $q(\cdot|x,a)$ is conservative and stable, i.e.,
\begin{eqnarray}\label{ZyExponential57}
q(S|x,a)=0,~\bar{q}_x=\sup_{a\in A}q_x(a)<\infty,
\end{eqnarray}
where $q_x(a):=-q(\{x\}|x,a).$ The signed kernel $q$ is often called the
transition rate. For simplicity and to fix ideas, we do
not consider the case of different admissible action spaces at
different states. Practically, the case of state-dependent admissible action spaces can be often reduced to the current setup by assigning a cost rate of $\infty$ at an inadmissible action, c.f. p.402 of \cite{Feinberg:2005}.

Let us take the sample space $\Omega$ by adjoining to the
countable product space $S\times((0,\infty)\times S)^\infty$ the
sequences of the form
$(x_0,\theta_1,\dots,\theta_n,x_n,\infty,x_\infty,\infty,x_\infty,\dots),$
where $x_0,x_1,\dots,x_n$ belong to $S$,
$\theta_1,\dots,\theta_n$ belong to $(0,\infty),$ and
$x_{\infty}\notin S$ is the isolated point. We equip $\Omega$ with
its Borel $\sigma$-algebra $\cal F$.

Let $t_0(\omega):=0=:\theta_0,$ and for each $n\geq 0$, and each
element $\omega:=(x_0,\theta_1,x_1,\theta_2,\dots)\in \Omega$, let
\begin{eqnarray*}
t_n(\omega)&:=&t_{n-1}(\omega)+\theta_n,
\end{eqnarray*}
and
\begin{eqnarray*}
t_\infty(\omega):=\lim_{n\rightarrow\infty}t_n(\omega).
\end{eqnarray*}
Obviously, $(t_n(\omega))$ are measurable mappings on $(\Omega,{\cal
F})$. In what follows, we often omit the argument $\omega\in
\Omega$ from the presentation for simplicity. Also, we regard
$x_n$ and $\theta_{n+1}$ as the coordinate variables, and note
that the pairs $\{t_n,x_n\}$ form a marked point process with the
internal history $\{{\cal F}_t\}_{t\ge 0},$ i.e., the filtration
generated by $\{t_n,x_n\}$; see Chapter 4 of \cite{Kitaev:1995}
for greater details. The marked point process $\{t_n,x_n\}$
defines the stochastic process on $(\Omega,{\cal F})$ of interest
$\{\xi_t,t\ge 0\}$ by
\begin{eqnarray}\label{ZhangExponentialGZCTMDPdefxit}
\xi_t=\sum_{n\ge 0}I\{t_n\le t<t_{n+1}\}x_n+I\{t_\infty\le
t\}x_\infty.
\end{eqnarray}
Here we accept $0\cdot x:=0$ and $1\cdot x:=x$ for each $x\in S_\infty,$ and below we denote
$S_{\infty}:=S\bigcup\{x_\infty\}$.

\begin{definition}
A (history-dependent) policy $\pi$ for the CTMDP is given by a
sequence $(\pi_n)$ such that, for each $n=0,1,2,\dots,$
$\pi_n(da|x_0,\theta_1,\dots,x_{n},s)$ is a stochastic kernel on
$A$, and for each $\omega=(x_0,\theta_1,x_1,\theta_2,\dots)\in
\Omega$, $t> 0,$
\begin{eqnarray*}
\pi(da|\omega,t)&=&I\{t\ge t_\infty\}\delta_{a_\infty}(da)+
\sum_{n=0}^\infty I\{t_n< t\le
t_{n+1}\}\pi_{n}(da|x_0,\theta_1,\dots,\theta_n,x_n, t-t_n),
\end{eqnarray*}
where $a_\infty\notin A$ is some isolated point. A policy
$\pi=(\pi_n)$  is called Markov if, with slight abuse of
notations, each of the stochastic kernels $\pi_n$ reads $
\pi_n(da|x_0,\theta_1,\dots,x_{n},s)=\pi_n(da|x_{n},s).$ A
Markov policy is further called deterministic if the stochastic
kernels $\pi_n(da|x_{n},s)$ all degenerate. A policy
$\pi=(\pi_n)$  is called stationary if, with slight abuse of
notations, each of the stochastic kernels $\pi_n$ reads $
\pi_n(da|x_0,\theta_1,\dots,x_{n},s)=\pi(da|x_{n}). $ A
stationary policy is further called deterministic if $
\pi_n(da|x_0,\theta_1,\dots,x_{n},s)=\delta_{\{f(x_{n})\}}(da) $
for some measurable mapping $f$ from $S$ to $A$. We shall identify such a deterministic stationary policy by the underlying measurable mapping $f$.
\end{definition}
The class of all policies for the CTMDP model is denoted by $\Pi.$

Under a policy $\pi:=(\pi_n)\in \Pi$, we define the
following random measure on $S\times (0,\infty)$
\begin{eqnarray*}
\nu^\pi(dy, dt)&:=& \int_A \tilde{q}(dy|\xi_{t-}(\omega),a)\pi(da|\omega,t)dt\nonumber\\
&=&\sum_{n\ge 0}\int_A
\tilde{q}(dy|x_n,a)\pi_{n}(da|x_0,\theta_1,\dots, \theta_n,
x_n,t-t_n)I\{t_n< t\le t_{n+1}\}dt
\end{eqnarray*}
with $q_{x_\infty}(a_\infty)=q(dy|x_\infty,a_\infty):=0=:q_{x_\infty}(a)$ for each $a\in A.$ Then for each given initial distribution $\gamma$ on ${\cal B}(S),$ there exists a unique probability measure ${P}^\pi_\gamma$ such
that
\begin{eqnarray*}
{P}_{\gamma}^\pi(x_0\in dx)=\gamma(dx),
\end{eqnarray*}
and with
respect to $P_\gamma^\pi,$ $\nu^\pi$ is the dual predictable
projection of the random measure associated with the marked point
process $\{t_n,x_n\}$; see \cite{Jacod:1975,Kitaev:1995}. The
process $\{\xi_t\}$ defined by (\ref{ZhangExponentialGZCTMDPdefxit}) under the
probability measure ${}{P}_\gamma^\pi$ is called a CTMDP. Below,
when $\gamma$ is a Dirac measure concentrated at $x\in S,$
we use the denotation ${}{P}_x^\pi.$ Expectations with respect to
${}{P}_\gamma^\pi$ and ${}{P}_x^\pi$ are denoted as
${}{E}_{\gamma}^\pi$ and ${}{E}_{x}^\pi,$ respectively.

The following remark follows from \cite{Jacod:1975}.
\begin{remark}\label{ZyExponentialRemark1}
Under a fixed policy $\pi=(\pi_n)$, the conditional distribution of $(\theta_{n+1},x_{n+1})$ with the condition on $x_0,\theta_1,x_1,\dots,\theta_{n},x_n$ is given on $\{\omega:x_n(\omega)\in S\}$ by
\begin{eqnarray*}
&&P_\gamma^\pi(\theta_{n+1}\in \Gamma_1,~x_{n+1}\in \Gamma_2|x_0,\theta_1,x_1,\dots,\theta_{n},x_n)\\
&=&\int_{\Gamma_1}e^{-\int_0^t \int_A q_{x_n}(a)\pi_n(da|x_0,\theta_1,\dots,\theta_n,x_n,s)ds}\int_{A}\tilde{q}(\Gamma_2|x_n,a)\pi_n(da|x_0,\theta_1,\dots,\theta_n,x_n,t)dt,\\
&&~\forall~\Gamma_1\in{\cal B}((0,\infty)),~\Gamma_2\in{\cal B}(S);\\
&&P_\gamma^\pi(\theta_{n+1}=\infty,~x_{n+1}=x_\infty|x_0,\theta_1,x_1,\dots,\theta_{n},x_n)=e^{-\int_0^\infty  \int_A q_{x_n}(a)\pi_n(da|x_0,\theta_1,\dots,\theta_n,x_n,s)ds},
\end{eqnarray*}
and given on $\{\omega:x_n(\omega)=x_\infty\}$ by
\begin{eqnarray*}
P_\gamma^\pi(\theta_{n+1}=\infty,~x_{n+1}=x_\infty|x_0,\theta_1,x_1,\dots,\theta_{n},x_n)=1.
\end{eqnarray*}
\end{remark}

Let the cost rate be given by a $[0,\infty)$-valued measurable function $c$ on $S\times A.$ In this paper, we study the following optimal control problem:
\begin{eqnarray}\label{ZyExponChap1}
\mbox{Minimize over $\pi\in \Pi$:}&&
E_x^\pi\left[e^{\int_0^\infty \int_A
c_0(\xi_t,a)\pi(da|\omega,t)dt}\right]=:V(x,\pi),~x\in S.
\end{eqnarray}
Here and below, we put $c(x_\infty,a):=0$ for each $a\in A.$

The CTMDP problem (\ref{ZyExponChap1}) is equivalent to minimizing the certainty equivalent of the total cost with respect to the exponential utility for a risk averse decision maker, see \cite{Cavazos:2000,Ghosh:2014,Howard:1972}.

In what follows, we refer the CTMDP problem (\ref{ZyExponChap1}) with the exponential utility to as the CTMDP model $\{S,A,q,c\}.$

\begin{definition}
A policy $\pi^\ast$ is called optimal for problem (\ref{ZyExponChap1}) if
\begin{eqnarray*}
V(x,\pi^\ast)=\inf_{\pi\in \Pi}V(x,\pi)=:V^\ast(x),~\forall~x\in S.
\end{eqnarray*}
\end{definition}
Evidently, $V^\ast(x)\ge 1$ for each $x\in S.$

One powerful method of reducing a CTMDP to a DTMDP is due to Yushkevich \cite{Yushkevich:1980}, which considers the case of a linear utility. However, the induced DTMDP is with a more complicated action space, so that a deterministic stationary strategy in the induced DTMDP in general does not give a deterministic stationary policy for the CTMDP model, but gives a specific Markov policy. This approach has been further developed to study piecewise deterministic Markov decision processes, see the books \cite{BauerleRieder:2011,Costa:2013,Davis:1993}. In Section \ref{ZyExponentialSec02}, as a preparation for our main optimality result, we shall develop this method for the case of exponential utility. In contrast to the linear utility case, now both the state and action spaces of the reduced DTMDP are more complicated than those of the CTMDP; e.g., a deterministic stationary strategy in this reduced DTMDP will no longer give a Markov policy for the CTMDP model. A further reduction to a simpler DTMDP with the same state and action spaces as the original CTMDP will be given in a subsequent section.

\section{Discrete-time Markov Decision Process with exponential utility}\label{SICONDTMDP}

To serve the investigations of the CTMDP, in this section we present briefly the dynamic programming approach for the DTMDP model (with exponential utility). The presented results are mostly related to \cite{Jaskiewicz:2008}, which however, is based on the compactness-continuity condition. For our purpose, we would not assume any compactness-continuity condition here, except for Proposition \ref{ZyExponentialProposition02}, and would need to consider a slightly more general cost function, as compared to \cite{Jaskiewicz:2008}. Without assuming the compactness-continuity condition, the dynamic programming approach for the DTMDP model was partially studied in \cite{Bertsekas:1978}, which dealt with a bounded cost function and mainly a finite horizon, see p.90 and Section 11.3 therein. Instead of the dynamic programming approach,  \cite{Fainberg:1982,Feinberg:1996} developed a different method for studying a rather general class of DTMDP problems in Borel spaces. That method, which can be traced back to Girsanov, is based on the investigations of strategic measures, and does not give all the results we would need here.

Consider a discrete-time Markov decision process with the following primitives:
\begin{itemize}
\item $\textbf{X}$ is a nonempty Borel state space.
\item $\textbf{A}$ is a nonempty Borel action space.
\item $p(dy|x,a)$ is a stochastic kernel on ${\cal B}(\textbf{X})$ given $(x,a)\in \textbf{X}\times\textbf{A}$.
\item $l$ a $[0,\infty]$-valued measurable cost function on $\textbf{X}\times\textbf{A}\times\textbf{X}.$
\end{itemize}
Let us denote for each $n=1,2,\dots,\infty,$ $\textbf{H}_n:=\textbf{X}\times(\textbf{A}\times\textbf{X})^n$ and $\textbf{H}_{0}:=\textbf{X}.$
A strategy $\sigma=(\sigma_n)_{n=0}^\infty$ in the DTMDP is given by a sequence of stochastic kernels $\sigma_n(da|h_{n})$ on ${\cal B}(\textbf{A})$ from $h_{n}\in \textbf{H}_{n}$ for $n=0,1,2,\dots.$ A strategy $\sigma=(\sigma_n)$ is called deterministic Markov if for each $n=0,1,2,\dots,$ $\sigma_n(da|h_{n})=\delta_{\{\varphi_n(x_{n})\}}(da)$, where $\varphi_{n}$ is an $\textbf{A}$-valued measurable mapping on $\textbf{X}.$ We identify such a deterministic Markov strategy with $(\varphi_n).$

Let $\Sigma$ be the space of strategies, and $\Sigma_{DM}$ be the space of all deterministic strategies for the DTMDP.

Let the controlled and controlling process be denoted by $\{Y_n, n=0,1,\dots,\infty\}$ and $\{A_n,n=0,1,\dots,\infty\}$. Here, for each $n=0,1,\dots,$ $Y_n$ is the projection of $\textbf{H}_\infty$ to the $2n+1$st coordinate, and $A_n$ to the $2n+2$nd coordinate.

Under a strategy $\sigma=(\sigma_n)$ and a given initial probability distribution $\nu$ on $(\textbf{X},{\cal B}(\textbf{X}))$, by the Ionescu-Tulcea theorem, c.f., \cite{Hernandez-Lerma:1996,Piunovskiy:1997}, one can construct a probability measure $\textbf{P}_\nu^\sigma$ on $(\textbf{H}_\infty,{\cal B}(\textbf{H}_\infty))$ such that
\begin{eqnarray*}
&&\textbf{P}_\nu^\sigma(Y_0\in dx)=\nu(dx),\\
&&\textbf{P}_\nu^\sigma(A_n\in da|Y_0,A_0,\dots,Y_n)=\sigma_{n}(da|Y_0,A_0,\dots,Y_n),~n=0,1,\dots,\\
&&\textbf{P}_\nu^\sigma(Y_{n+1}\in dx|Y_0,A_0,\dots,Y_n,A_n)=p(dx|Y_n,A_n),~n=0,1,\dots.
\end{eqnarray*}
As usual, equalities involving conditional expectations and probabilities are understood in the almost sure sense.
The probability measure $\textbf{P}_\nu^\sigma$ is called a strategic measure for the DTMDP. The expectation taken with respect to $\textbf{P}_\nu^\sigma$ is denoted by $\textbf{E}_\nu^\sigma.$ When $\nu$ is concentrated on the singleton $\{x\}$, $\textbf{P}_\nu^\sigma$ and $\textbf{E}_\nu^\sigma$ are written as $\textbf{P}_x^\sigma$ and $\textbf{E}_x^\sigma.$

Consider the optimal control problem
\begin{eqnarray}\label{ZyExponentialProblem2}
\mbox{Minimize over $\sigma$}:&& \textbf{E}_x^\sigma\left[e^{\sum_{n=0}^\infty l(Y_n,A_n,Y_{n+1})}\right]=:\textbf{V}(x,\sigma),~x\in \textbf{X}.
\end{eqnarray}
We denote the value function of problem (\ref{ZyExponentialProblem2}) by $\textbf{V}^\ast$. Then a strategy $\sigma^\ast$ is called optimal for problem  (\ref{ZyExponentialProblem2}) if $\textbf{V}(x,\sigma^\ast)=\textbf{V}^\ast(x)$ for each $x\in \textbf{X}.$
We refer problem (\ref{ZyExponentialProblem2}) to as the DTMDP model $\{\textbf{X},\textbf{A},p,l\}$ (with the exponential utility).

Note that
\begin{eqnarray}\label{ZyExponential08}
\textbf{V}^\ast(x)\ge 1,~\forall~x\in\textbf{X}.
\end{eqnarray}

One can write
\begin{eqnarray}\label{ZyExponential05}
\textbf{E}_x^\sigma\left[e^{\sum_{n=0}^\infty l(Y_n,A_n,Y_{n+1})}\right]=\int_{\textbf{H}_\infty}e^{\sum_{n=0}^\infty l(Y_n,A_n,Y_{n+1})}\textbf{P}_x^\sigma(dh),~x\in \textbf{X}.
\end{eqnarray}
Then $\textbf{V}(x,\cdot)$ is a measurable criterion in the sense of \cite{Fainberg:1982,Feinberg:1996}.
In view of this, that $\textbf{V}^\ast$ is a lower semianalytic function on $\textbf{X}$ immediately follows from Theorem 4.2 of \cite{Feinberg:1996}.

\begin{proposition}\label{ZyExponentialTheorem02}
The function $\textbf{V}^\ast$ is an $[1,\infty]$-valued lower semianalytic solution to
\begin{eqnarray}\label{ZyExponential02}
\textbf{V}(x)=\inf_{a\in \textbf{A}}\left\{\int_{\textbf{X}}p(dy|x,a)e^{l(x,a,y)}\textbf{V}(y)\right\},~x\in \textbf{X}.
\end{eqnarray}
\end{proposition}
\par\noindent\textit{Proof.}
Let us verify that $\textbf{V}^\ast$ solves the optimality equation (\ref{ZyExponential02}) as follows. (The other assertions in the statement of this theorem have been justified in the previous discussions.)

Let $\sigma$ be an arbitrarily fixed strategy. For each $x\in \textbf{X}$ and $b\in\textbf{A}$, consider the shifted strategy $\sigma^{(x,b)}=(\sigma^{(x,b)}_n)$ by $\sigma^{(x,b)}_n(da|h_n)=\sigma_{n+1}(da|x,b,h_n)$ for each $n=0,1,\dots,$ and $h_n\in \textbf{H}_n.$ Then for each fixed $x\in \textbf{X}$,
\begin{eqnarray*}
\textbf{V}^\ast(x)\le \textbf{V}(x,\sigma^{(x,b)})=\left.\textbf{E}_x^\sigma \left[e^{\sum_{n=1}^\infty l(Y_n,A_n,Y_{n+1})}\right|Y_0=x,A_0=b,Y_1=y \right],
 \end{eqnarray*}
where the equality holds almost surely with respect to $p(dy|x,b)\sigma_0(db|x)$. Now for each $x\in \textbf{X}$,
\begin{eqnarray*}
\textbf{V}(x,\sigma)&=&\textbf{E}_x^\sigma\left[e^{l(Y_0,A_0,Y_1)}e^{\sum_{n=1}^\infty l(Y_n,A_n,Y_{n+1})} \right]\\
&=&\left.\textbf{E}_x^\sigma\left[e^{l(Y_0,A_0,Y_1)}\textbf{E}_x^\sigma\left[e^{\sum_{n=1}^\infty l(Y_n,A_n,Y_{n+1})}\right|Y_0,A_0,Y_1 \right]\right]\\
&=&\left.\int_{\textbf{A}}\int_{\textbf{X}} p(dy|x,b)e^{l(x,b,y)}\textbf{E}_x^\sigma \left[e^{\sum_{n=1}^\infty l(Y_n,A_n,Y_{n+1})}\right|Y_0=x,A_0=b,Y_1=y \right]\sigma_0(db|x)\\
&\ge &\int_{\textbf{A}} \int_{\textbf{X}} p(dy|x,b)e^{l(x,b,y)}\textbf{V}^\ast(y)\sigma_0(db|x)\\
&\ge& \inf_{a\in \textbf{A}}\left\{\int_{\textbf{X}} p(dy|x,a)e^{l(x,a,y)}\textbf{V}^\ast(y)\right\},
\end{eqnarray*}
where and below integrals such as those in the above inequalities are well defined because $\textbf{V}^\ast$ is lower semianalytic, see Lemma 7.30 as well as Proposition 7.48 of \cite{Bertsekas:1978}.
Thus,
\begin{eqnarray}\label{ZyExponential01}
\textbf{V}^\ast(x)\ge \inf_{a\in \textbf{A}}\left\{\int_{\textbf{X}} p(dy|x,a)e^{l(x,a,y)}\textbf{V}^\ast(y)\right\},~x\in \textbf{X}.
\end{eqnarray}

Next, we establish the opposite direction of the above inequality.
Let $z\in \textbf{X}$, $b\in \textbf{A}$ and $\epsilon>0$  be arbitrarily fixed. It follows from Theorem 3.1 of \cite{Fainberg:1982} (c.f. Chapter 3 of \cite{Dynkin:1979}) that there exists a strategy $\hat{\sigma}$ such that
\begin{eqnarray*}
\ln \textbf{V}^\ast(x)+\epsilon\ge \ln \textbf{E}_x^{\hat{\sigma}}\left[e^{\sum_{n=0}^\infty l(Y_n,A_n,Y_{n+1})}\right]=\ln \textbf{V}(x,\hat{\sigma})
\end{eqnarray*}
for almost all $x\in \textbf{X}$ with respect to $p(dx|z,b).$ Consider the strategy $\sigma'=(\delta_{\{b\}},\hat{\sigma}_0,\hat{\sigma}_1,\dots)$. Then
\begin{eqnarray*}
\textbf{V}^\ast(z)&\le& \textbf{E}_z^{\sigma'}\left[e^{\sum_{n=0}^\infty l(Y_n,A_n,Y_{n+1})}\right]=\int_{\textbf{X}}p(dy|z,b)e^{l(z,b,y)}e^{\ln \textbf{V}(y,\hat{\sigma})}\\
&\le& \int_{\textbf{X}}p(dy|z,b)e^{l(z,b,y)} e^{\ln \textbf{V}^\ast(y)+\epsilon}\\
&=&e^\epsilon\int_{\textbf{X}}p(dy|z,b)e^{l(z,b,y)}\textbf{V}^\ast(y).
\end{eqnarray*}
Since $\epsilon>0$, $z\in \textbf{X}$ and $b\in \textbf{A}$ are arbitrarily fixed, it follows from the above that
\begin{eqnarray*}
\textbf{V}^\ast(x)\le \inf_{a\in \textbf{A}}\left\{\int_{\textbf{X}} p(dy|x,a)e^{l(x,a,y)}\textbf{V}^\ast(y)\right\},~x\in \textbf{X}.
\end{eqnarray*}
Combing this with (\ref{ZyExponential01}), we see that the statement holds.
$\hfill\Box$
\bigskip

\begin{proposition}\label{ZyExponentialTheorem01}
For each $x\in\textbf{X},$
\begin{eqnarray*}
\inf_{\sigma\in \Sigma_{DM}}\textbf{V}(x,\sigma)=\inf_{\sigma\in \Sigma}\textbf{V}(x,\sigma),
\end{eqnarray*}
where we recall that $\Sigma_{DM}$ is the set of all deterministic Markov strategies for the DTMDP.
\end{proposition}

The proof of this proposition is based on the next result.
\begin{lemma}\label{ZyExponentialLemma01L}
Let $f$ be a $[1,\infty]$-valued lower semianalytic function on the Borel space $\textbf{X}\times\textbf{A},$ and $f^\ast$ be a function on $\textbf{X}$ defined by $f^\ast(x)=\inf_{a\in \textbf{A}}f(x,a)$ for each $x\in \textbf{X}.$ Then for each $\epsilon>0,$ there exists an analytically measurable mapping $\varphi$ from $\textbf{X}$ to $\textbf{A}$ such that
\begin{eqnarray*}
f(x,\varphi(x))\le f^\ast(x) e^{\epsilon},~\forall~x\in \textbf{X}.
\end{eqnarray*}
\end{lemma}
\par\noindent\textit{Proof.} The reasoning of Proposition 7.50 of \cite{Bertsekas:1978} can be easily modified to prove the statement of this lemma. The details are omitted. $\hfill\Box$
\bigskip

Now we are in position to prove Proposition \ref{ZyExponentialTheorem01}.
\bigskip

\par\noindent\textit{Proof of Proposition \ref{ZyExponentialTheorem01}.} Let $\epsilon>0$ and $x_0\in \textbf{X}$ be arbitrarily fixed. We first show that there exists a deterministic Markov strategy $\sigma$ such that
\begin{eqnarray}\label{ZyExponential09}
\textbf{V}(x_0,\sigma)\le \textbf{V}^\ast(x_0)e^\epsilon
\end{eqnarray}
as follows.

Let $(\epsilon_k)$ be a sequence of positive constants such that $\sum_{k=0}^\infty \epsilon_k=\epsilon.$ By Proposition \ref{ZyExponentialTheorem02}, there is a Borel measurable mapping $\varphi_0$ from $\textbf{X}$ to $\textbf{A}$ such that
\begin{eqnarray}\label{ZyExponential06}
\textbf{V}^\ast(x_0)\ge e^{-\epsilon_0}\int_{\textbf{X}}p(dx_1|x_0,\varphi_0(x_0))e^{l(x_0,\varphi_0(x_0),x_1)}\textbf{V}^\ast(x_1).
\end{eqnarray}
(Remember, the above inequality is only required to hold for the fixed $x_0\in \textbf{X}.$)
By Lemma \ref{ZyExponentialLemma01L} and Proposition \ref{ZyExponentialTheorem02}, for each $k=1,2,\dots,$ there exists an analytically measurable mapping $\tilde{\varphi}_k$ from $\textbf{X}$ to $\textbf{A}$ such that
\begin{eqnarray}\label{ZyExponential07}
\textbf{V}^\ast(x)\ge e^{-\epsilon_k}\int_{\textbf{X}}p(dy|x,\tilde{\varphi}_k(x))e^{l(x,\tilde{\varphi}_k(x),y)}\textbf{V}^\ast(y),~\forall~x\in \textbf{X}.
\end{eqnarray}
Let $\varphi_1$ be the Borel measurable modification of $\tilde{\varphi}_1$ with respect to the probability measure $p(\cdot|x_0,\varphi_0(x_0))$. Then
\begin{eqnarray*}
\textbf{V}^\ast(x_0)&\ge&  e^{-\epsilon_0-\epsilon_1}\int_{\textbf{X}}p(dx_1|x_0,\varphi_0(x_0))e^{l(x_0,\varphi_0(x_0),x_1)}\int_{\textbf{X}}p(dx_2|x_1,\tilde{\varphi}_1(x_1))e^{l(x_1,\tilde{\varphi}_1(x_1),x_2)}\textbf{V}^\ast(x_2)\\
&\ge &e^{-\epsilon_0-\epsilon_1}\int_{\textbf{X}}p(dx_1|x_0,\varphi_0(x_0))e^{l(x_0,\varphi_0(x_0),x_1)}\int_{\textbf{X}}p(dx_2|x_1,\varphi_1(x_1))e^{l(x_1,\varphi_1(x_1),x_2)},
\end{eqnarray*}
where the first inequality is by (\ref{ZyExponential06}) and (\ref{ZyExponential07}), and the second inequality is by (\ref{ZyExponential08}).
Inductively, for each $k=2,3,\dots,$ let $\varphi_k$ be a Borel measurable modification of $\tilde{\varphi}_k$ with respect to the probability measure $\int_{\textbf{X}^{k-1}}p(\cdot|x_{k-1},\varphi_{k-1}(x_{k-1}))p(dx_{k-1}|x_{k-2},\varphi_{k-2}(x_{k-2}))\dots p(dx_1|x_0,\varphi_0(x_0)).$ Let $\sigma=(\varphi_n)_{n=0}^\infty$ be the deterministic Markov strategy. Then
\begin{eqnarray*}
\textbf{V}^\ast(x_0)&\ge&  e^{-\sum_{n=0}^k\epsilon_n} \textbf{E}_{x_0}^\sigma\left[e^{\sum_{n=0}^k l(Y_n,A_n,Y_{n+1})}\right],~\forall~k=0,1,\dots.
\end{eqnarray*}
By passing to the limit as $k\rightarrow\infty$ in the above inequality, we see that (\ref{ZyExponential09}) holds.

Now the statement of the theorem follows from (\ref{ZyExponential09}) and the arbitrariness of $\epsilon>0$ and $x_0\in \textbf{X}$.
  $\hfill\Box$
\bigskip

\begin{proposition}\label{ZyExponentialProposition01} The following two assertions hold.
\begin{itemize}
\item[(a)] Let $\textbf{U}$ be a $[1,\infty]$-valued lower semianalytic function on $\textbf{X}$. If
\begin{eqnarray*}
\textbf{U}(x)\ge \inf_{a\in \textbf{A}}\left\{\int_{\textbf{X}}p(dy|x,a)e^{l(x,a,y)}\textbf{U}(y)\right\},~\forall~x\in \textbf{X},
\end{eqnarray*}
then $\textbf{U}(x)\ge \textbf{V}^\ast(x)$ for each $x\in \textbf{X}.$ In particular, if the function $\textbf{U}$ satisfying the above relation is $[1,\infty)$-valued, then so is the value function $\textbf{V}^\ast.$
\item[(b)] Let $\varphi$ be a deterministic stationary strategy for the DTMDP model $\{\textbf{X},\textbf{A},p,l\}$. If
\begin{eqnarray}\label{ZyExponential025}
\textbf{V}^\ast(x)=\int_{\textbf{X}}p(dy|x,\varphi(x))e^{l(x,\varphi(x),y)}\textbf{V}^\ast(y),~\forall~x\in \textbf{X},
\end{eqnarray}
then $\textbf{V}^\ast(x)=\textbf{V}(x,\varphi)$ for each $x\in \textbf{X}.$
\end{itemize}
\end{proposition}
\par\noindent\textit{Proof.} (a) Let $x\in \textbf{X}$ and $\epsilon>0$ be fixed. Then by using Lemma \ref{ZyExponentialLemma01L}, one can follow the reasoning in the proof of Proposition \ref{ZyExponentialTheorem01}, and see the existence of a deterministic Markov strategy $\sigma$, which satisfies
\begin{eqnarray*}
\textbf{U}(x)\ge e^{-\epsilon} \textbf{V}(x,\sigma)\ge e^{-\epsilon}\textbf{V}^\ast(x).
\end{eqnarray*}
Since $\epsilon>0$ and $x\in \textbf{X}$ are arbitrarily fixed, the statement follows.

(b) Consider the given deterministic stationary strategy $\varphi.$ Let $x\in \textbf{X}$ be fixed. Then by simple iterations based on (\ref{ZyExponential025}), and keeping in mind that $\textbf{V}^\ast$ is $[1,\infty]$-valued, we see
\begin{eqnarray*}
\textbf{V}^\ast(x)= \textbf{E}_x^\varphi\left[e^{\sum_{n=0}^ml(Y_n,A_n,Y_{n+1})}\textbf{V}(Y_{m+1})\right]\ge \textbf{E}_x^\varphi\left[e^{\sum_{n=0}^ml(Y_n,A_n,Y_{n+1})}\right]
\end{eqnarray*}
for each $m=1,2,\dots.$ Thus, the statement holds after passing to the limit as $m\rightarrow\infty.$ $\hfill\Box$
\bigskip

\begin{condition}\label{ZyExponentialDTMDPCon1}
\begin{itemize}
\item[(a)] The function $l$ is lower semicontinuous on $\textbf{X}\times\textbf{A}\times\textbf{X}.$
\item[(b)] For each bounded continuous function $f$ on $\textbf{X}$, $\int_{\textbf{X}}f(y)p(dy|x,a)$ is continuous in $(x,a)\in\textbf{X}\times\textbf{A}.$
\item[(c)] The space $\textbf{A}$ is a compact Borel space.
\end{itemize}
\end{condition}

\begin{condition}\label{ZyExponentialDTMDPCon2}
\begin{itemize}
\item[(a)] The function $l(x,a,y)$ is lower semicontinuous in $a\in \textbf{A}$ for each $x,y\in \textbf{X}.$
\item[(b)] For each bounded measurable function $f$ on $\textbf{X}$ and each $x\in \textbf{X},$ $\int_{\textbf{X}}f(y)p(dy|x,a)$ is continuous in $a\in\textbf{A}.$
\item[(c)] The space $\textbf{A}$ is a compact Borel space.
\end{itemize}
\end{condition}

\begin{proposition}\label{ZyExponentialProposition02}
\begin{itemize}
\item[(a)] Suppose Condition \ref{ZyExponentialDTMDPCon1} is satisfied. Then the value function $\textbf{V}^\ast$ is the minimal $[1,\infty]$-valued lower semicontinuous solution to (\ref{ZyExponential02}).
\item[(b)] Suppose Condition \ref{ZyExponentialDTMDPCon2} is satisfied. Then the value function $\textbf{V}^\ast$ is the minimal $[1,\infty]$-valued measurable solution to (\ref{ZyExponential02}).
\item[(c)] Suppose Condition \ref{ZyExponentialDTMDPCon1} or Condition \ref{ZyExponentialDTMDPCon2} is satisfied. Let $\textbf{V}^{(0)}(x):=0$ for each $x\in \textbf{X}$, and for each $n=1,2,\dots,$
    \begin{eqnarray*}
    \textbf{V}^{(n)}(x):=\inf_{a\in A}\left\{\int_{\textbf{X}}p(dy|x,a)e^{l(x,a,y)}\textbf{V}^{(n-1)}(y)\right\},~\forall~x\in \textbf{X}.
    \end{eqnarray*}
Then $(\textbf{V}^{(n)}(x))$ increases to $\textbf{V}^\ast(x)$ for each $x\in \textbf{X}$, where $\textbf{V}^\ast$ is the value function for problem (\ref{ZyExponentialProblem2}). Furthermore, there exists a deterministic stationary strategy $\varphi$ satisfying (\ref{ZyExponential025}), and so in particular, there exists a deterministic stationary optimal strategy for the DTMDP problem (\ref{ZyExponentialProblem2}).
\end{itemize}
\end{proposition}
\par\noindent\textit{Proof.} This statement can be proved as in \cite{Jaskiewicz:2008}, even though in \cite{Jaskiewicz:2008} the function $l(x,a,y)$ does not depend on $y\in \textbf{X}$. $\hfill\Box$

\section{First reduction to a DTMDP model}\label{ZyExponentialSec02}
In this section, we reduce the risk-sensitive CTMDP to a risk-sensitive DTMDP with more complicated state and action spaces. As in \cite{Yushkevich:1980}, this is based on viewing a policy $\pi=(\pi_n)$ as a sequence of measurable mappings taking values in the (quotient) space of $\mathbb{P}(A)$-valued measurable mappings, where    and below, the set $\mathbb{P}(A)$ is the space of
probability measures on ${\cal B}(A)$, and is equipped with the
standard weak topology, so that $\mathbb{P}(A)$ is a Borel space, see Chapter 7 of \cite{Bertsekas:1978}. The details are as follows.

Let ${\cal R}$ denote the set of (Borel) measurable mappings
$\rho_t(da)$ from $t\in(0,\infty)\rightarrow \mathbb{P}(A).$ Here,
we do not distinguish between two measurable mappings in $t\in
(0,\infty)$ which coincide almost everywhere with respect to the
Lebesgue measure.

We endow ${\cal R}$ with the $\sigma$-algebra as the minimal one with respect to which, the function
\begin{eqnarray*}
\rho\in {\cal{R}}\rightarrow \int_0^\infty e^{-t}
f(t,\rho_t)dt
\end{eqnarray*}
is measurable in $\rho\in {\cal R}$ for each bounded measurable function $f$ on $(0,\infty)\times \mathbb{P}(A).$ Lemma 1 of \cite{Yushkevich:1980} asserts that ${\cal R}$ is a Borel space.

For the rest of this paper, it is convenient to introduce the
following notations. For each $\mu\in \mathbb{P}(A)$,
\begin{eqnarray*}
&& q_x(\mu):=\int_A q_x(a)\mu(da);\\
&&\tilde{q}(dy|x,\mu):= \int_A \tilde{q}(dy|x,a)\mu(da);\\
&&c(x,\mu):=\int_A c(x,a)\mu(da).
\end{eqnarray*}

It follows from \cite{Yushkevich:1980} that we can legitimately consider a DTMDP model $\{\textbf{X},\textbf{A},p,l\}$ with exponential utility with the following primitives, where all the functions and mappings are measurable.
\begin{itemize}
\item The state space is $\textbf{X}:=((0,\infty)\times S)\bigcup\{(\infty,x_\infty)\}$. Whenever the topology is concerned, $(\infty,x_\infty)$ is regarded as an isolated point in $\textbf{X}.$
\item The action space is $\textbf{A}:={\cal R}$.
\item The transition kernel $p$ on ${\cal B}(\textbf{X})$ from $\textbf{X}\times \textbf{A}$ is given for each $\rho\in \textbf{A}$ by
    \begin{eqnarray}\label{ZyExponential52}
    p(\Gamma_1\times\Gamma_2|(\theta,x),\rho)&:=&\int_{\Gamma_2} e^{-\int_0^t q_x(\rho_s)ds}\tilde{q}(\Gamma_1|x,\rho_t)dt,\nonumber\\
    &&~\forall~\Gamma_1\in {\cal B}(S),~\Gamma_2\in {\cal B}((0,\infty)),~x\in S,~\theta\in (0,\infty),\nonumber\\
     p(\{\infty\}\times\{x_\infty\}|(\theta,x),\rho)&:=&e^{-\int_0^\infty q_x(\rho_s)ds},~\forall~x\in S,~\theta\in (0,\infty);\nonumber\\
     p(\{(\infty,x_\infty)\}|(\infty,x_\infty),\rho)&:=&1.\nonumber\\
\end{eqnarray}
\item The cost function $l$ is a $[0,\infty]$-valued measurable function on $\textbf{X}\times\textbf{A}\times \textbf{X}$ given by
\begin{eqnarray}\label{ZyExponential55}
l((\theta,x),\rho,(\tau,y)):=\int_0^\infty I\{s<\tau\} c(x,\rho_s)ds=:\hat{l}(x,\rho,\tau),~\forall~((\theta,x),\rho,(\tau,y))\in \textbf{X}\times\textbf{A}\times \textbf{X}.
\end{eqnarray}
\end{itemize}
(Recall that $c(x_\infty,a):=0$ and $q_{x_\infty}(a):=0=:\tilde{q}(S|x_\infty,a)$ for each $a\in A.$) For each strategy $\sigma$ for the DTMDP $\{\textbf{X},\textbf{A},p,l\}$, the function $\textbf{V}(\cdot,\sigma)$ is defined by (\ref{ZyExponentialProblem2}).

The controlled process in the above DTMDP model  $\{\textbf{X},\textbf{A},p,l\}$  is denoted by $\{Y_n,n=0,1,\dots\}$, where $Y_n=(\Theta_n,X_n)$, and the controlling process is denoted by $\{A_n,n=0,1,\dots\}.$
Let $\Sigma_{DM}^0$ be the class of deterministic Markov strategies for the DTMDP model $\{\textbf{X},\textbf{A},p,l\}$ in the form $\sigma=(\varphi_n)$ where $\varphi_0((\theta,x))$ does not depend on $\theta\in (0,\infty)$ for each $x\in S.$

For each fixed $\hat{\theta}\in (0,\infty)$ and a deterministic Markov strategy $\sigma=(\varphi_n)$,
\begin{eqnarray}\label{ZyExponential010}
\textbf{V}((\hat{\theta},x),\sigma)=\textbf{V}((\hat{\theta},x),\hat{\sigma})=\textbf{V}((\theta,x),\hat{\sigma}),~\forall~x\in S,~\theta\in(0,\infty),
\end{eqnarray}
where $\hat{\sigma}=(\hat{\varphi}_0,\varphi_1,\varphi_2,\dots)$ with $\hat{\varphi}_0((\theta,x))=\varphi_0((\hat{\theta},x))$ for each $\theta\in (0,\infty).$ Indeed, since $\hat{\theta}$ is fixed, under the strategy $\hat{\sigma}$, the decision is made independently of the first coordinate of the initial state. This, together with the definitions of the transition kernel $p$ and the cost function $l$ given by (\ref{ZyExponential52}) and (\ref{ZyExponential55}), justifies (\ref{ZyExponential010}). See also Theorem 2 of \cite{Feinberg:2005}.

Proposition \ref{ZyExponentialTheorem01} and (\ref{ZyExponential010}) imply that
\begin{eqnarray*}
\textbf{V}^\ast((\theta,x))=\inf_{\sigma\in \Sigma_{DM}^0}\textbf{V}((\theta,x),\sigma),~\forall~x\in S,~\theta\in(0,\infty),
\end{eqnarray*}
where $\textbf{V}^\ast$ is the value function of the DTMDP problem (\ref{ZyExponentialProblem2}). This together with (\ref{ZyExponential010}) then leads to that
$\textbf{V}^\ast((\theta,x))$ does not depend on $\theta\in (0,\infty).$ Therefore, we write
$ \textbf{V}^\ast(x)$ instead of $\textbf{V}^\ast((\theta,x))$ and $\textbf{V}(x,\sigma)$ instead of $\textbf{V}((\theta,x),\sigma)$ when $\sigma$ is in $\Sigma_{DM}^0$. The previous equality now reads
\begin{eqnarray}\label{ZyExponential011}
\textbf{V}^\ast(x)=\inf_{\sigma\in \Sigma_{DM}^0}\textbf{V}(x,\sigma),~\forall~x\in S.
\end{eqnarray}

Consider a policy $\pi=(\pi_n)$ for the CTMDP model $\{S,A,q,c\}$. Note that each stochastic kernel
\begin{eqnarray*}
\pi_n(da|x_0,\theta_1,x_1,\theta_2,\dots,\theta_n,x_n,s)
\end{eqnarray*}
can be identified with a measurable mapping say $\pi_n(x_0,\theta_1,x_1,\theta_2,\dots,\theta_n,x_n)(s,da)$ from $S\times \textbf{X}^n$ to ${\cal R},$ and vice versa. Therefore, each policy $\pi=(\pi_n)$ for the CTMDP model $\{S,A,q,c\}$ is identified with a deterministic strategy denoted by $\sigma(\pi)$ for the DTMDP model $\{\textbf{X},\textbf{A},p,l\}$, where, under this strategy $\sigma(\pi)$, at the time step $n,$ the decision in $\textbf{A}$ is made based only on $X_0, Y_1,Y_2,\dots,Y_n$ and $n$, and is independent on $\Theta_0$ and the past actions. Therefore, under the policy $\pi=(\pi_n)$ for the CTMDP model $\{S,A,q,c\}$, for each $x\in S$ and $\theta\in(0,\infty)$,
\begin{eqnarray*}
V(x,\pi)&=&E_x^\pi\left[e^{\int_0^\infty \int_A c(x,a)\pi(da|\omega,t)dt}\right]=E_x^\pi\left[e^{\sum_{n=0}^\infty  \int_{0}^{\theta_{n+1}}\int_A c(x_n,a)\pi_n(da|x_0,\theta_1,\dots,\theta_n,x_n,s)ds}\right]\\
&=&\textbf{E}_{(\theta,x)}^{\sigma(\pi)}\left[e^{\sum_{n=0}^\infty l(Y_n,A_n,Y_{n+1})}\right]=\textbf{E}_{(\theta,x)}^{\sigma(\pi)}\left[e^{\sum_{n=0}^\infty \hat{l}(X_n,A_n,\Theta_{n+1})}\right]=\textbf{V}(x,\sigma(\pi)),
\end{eqnarray*}
c.f. Remark \ref{ZyExponentialRemark1}, (\ref{ZyExponential52}) and (\ref{ZyExponential55}) for the third equality, and Theorem 2 of \cite{Feinberg:2005} for the last equality.
Thus,
\begin{eqnarray*}
V^\ast(x)\ge \textbf{V}^\ast(x),~\forall~x\in S.
\end{eqnarray*}

On the other hand, each deterministic Markov strategy $\sigma=(\varphi_n)\in \Sigma_{DM}^0$ can be identified with a policy say $\pi(\sigma)$ for the CTMDP model $\{S,A,q,c\}$ such that
\begin{eqnarray*}
V(x,\pi(\sigma))=\textbf{V}(x,\sigma),~\forall~x\in S.
\end{eqnarray*}
This and (\ref{ZyExponential011}) imply
\begin{eqnarray*}
V^\ast(x)\le \textbf{V}^\ast(x),~\forall~x\in S.
\end{eqnarray*}

Now we come to the following conclusion.
\begin{theorem}\label{ZyExponentialTheorem05}
The value function $V^\ast$ for the CTMDP problem (\ref{ZyExponChap1}) is lower semianalytic on $S$ and satisfies
\begin{eqnarray*}
V^\ast(x)=\textbf{V}^\ast(x),~\forall~x\in S.
\end{eqnarray*}
Furthermore, $V^\ast$ satisfies
\begin{eqnarray}\label{ZyExponential015}
V^\ast(x)
&=&\inf_{\rho\in {\cal R}}\left\{\int_0^\infty e^{-\int_0^\tau (q_x(\rho_s)-c(x,\rho_s))ds} \left(\int_S V^\ast(y)\tilde{q}(dy|x,\rho_\tau)\right)d\tau +e^{-\int_0^\infty q_x(\rho_s)ds}e^{\int_0^\infty c(x,\rho_s)ds} \right\},\nonumber\\
&&\forall~x\in S.
\end{eqnarray}
\end{theorem}
\par\noindent\textit{Proof.} The equality between $V^\ast$ and $\textbf{V}^\ast$ on $S$ follows from the discussions above the theorem. The other parts of the statement are then by Proposition \ref{ZyExponentialTheorem02}. $\hfill\Box$\bigskip

In (\ref{ZyExponential015}), it is possible that $\int_0^\infty q_x(\rho_s)ds=\infty=\int_0^\infty c(x,\rho_s)ds,$ so that by  (\ref{ZyExponential56})
\begin{eqnarray*}
e^{-\int_0^\infty q_x(\rho_s)ds}e^{\int_0^\infty c(x,\rho_s)ds}=0\cdot\infty=0\ne \infty=e^{\infty}=e^{-\infty+\infty}=e^{-\int_0^\infty q_x(\rho_s)ds+\int_0^\infty c(x,\rho_s)ds}.
\end{eqnarray*}
On the other hand, by (\ref{ZyExponential57}),
\begin{eqnarray*}
e^{-\int_0^\tau (q_x(\rho_s)-c(x,\rho_s))ds}=e^{-\int_0^\tau q_x(\rho_s)ds}e^{\int_0^\tau c(x,\rho_s)ds}
\end{eqnarray*}
for each $\tau\in(0,\infty)$.

\section{Optimality result}\label{ZyExponentialSect03}
In this section, we establish the optimality equation for the CTMDP problem (\ref{ZyExponChap1}); show, under some compactness-continuity conditions, the existence of a deterministic stationary optimal policy for problem (\ref{ZyExponChap1}); and further reduce the CTMDP model $\{S,A,q,c\}$ to a simpler DTMDP model with the same state and action spaces as the CTMDP, in contrast to the DTMDP model in Section \ref{ZyExponentialSec02}. As a corollary of this reduction, we formulate the value iteration algorithm for problem (\ref{ZyExponChap1}).

\subsection{Optimality equation}
In this subsection, we establish the optimality equation satisfied by the value function $V^\ast$ of the CTMDP problem (\ref{ZyExponChap1}). This is done based on more detailed investigations of (\ref{ZyExponential015}).
\begin{theorem}\label{ZyExponentialTheorem06}
\begin{itemize}
\item[(a)] The value function $V^\ast$ of the CTMDP problem (\ref{ZyExponChap1}) is an $[1,\infty]$-valued lower semianalytic function satisfying
\begin{eqnarray}\label{ZyExponential036}
0=\inf_{a\in A}\left\{c(x,a)V^\ast(x)+\int_S q(dy|x,a)V^\ast(y)\right\}
\end{eqnarray}
for each $x\in S$ such that $V^\ast(x)<\infty.$
\item[(b)] If a deterministic stationary policy $\varphi$ for the CTMDP model $\{S,A,q,c\}$ satisfies
\begin{eqnarray}\label{ZyExponential026}
0=\inf_{a\in A}\left\{c(x,a)V^\ast(x)+\int_S q(dy|x,a)V^\ast(y)\right\}=c(x,\varphi(x))V^\ast(x)+\int_S q(dy|x,\varphi(x))V^\ast(y)
\end{eqnarray}
for each $x\in S$ such that $V^\ast(x)<\infty,$ then the deterministic stationary policy $\varphi$ is optimal for the CTMDP problem (\ref{ZyExponChap1}). (The definition of $\varphi$ on $\{x\in S:~V^\ast(x)=\infty\}$ can be put arbitrarily, so long $\varphi$ is measurable on $S$.)
\end{itemize}
\end{theorem}
We call (\ref{ZyExponential036}) the optimality equation for the CTMDP problem (\ref{ZyExponChap1}), and
call an $[1,\infty]$-valued lower semianalytic function $V$ a solution to the optimality equation (\ref{ZyExponential036}) if it satisfies (\ref{ZyExponential036}) with $V^\ast$ being replaced by $V$ for each $x\in S$, where $V(x)<\infty.$
To guarantee the existence of such a deterministic stationary policy $\varphi$ as in Theorem \ref{ZyExponentialTheorem06}(b), in the next subsection, we shall impose some compactness-continuity conditions, under which the value function will be seen to be measurable or lower semicontinuous.

We postpone the proof of Theorem \ref{ZyExponentialTheorem06} after several lemmas.
\begin{lemma}\label{ZyExponentialLemma02L}
For each $x\in S$ and $\rho\in {\cal R}$,
\begin{eqnarray*}
t\in[0,\infty)\rightarrow \int_0^t e^{-\int_0^\tau (q_x(\rho_s)-c(x,\rho_s))ds}\int_{S}V^\ast(y)\tilde{q}(dy|x,\rho_\tau)d\tau +e^{-\int_0^t (q_x(\rho_s)-c(x,\rho_s))ds}V^\ast(x)
\end{eqnarray*}
is monotone nondecreasing in $t\in[0,\infty)$.
\end{lemma}
\par\noindent\textit{Proof.} Let $0\le t_1<t_2<\infty$ be arbitrarily fixed. For the statement of the lemma, it suffices to show
\begin{eqnarray}\label{ZyExponential012}
&&\int_0^{t_2} e^{-\int_0^\tau (q_x(\rho_s)-c(x,\rho_s))ds}\int_{S}V^\ast(y)\tilde{q}(dy|x,\rho_\tau)d\tau +e^{-\int_0^{t_2} (q_x(\rho_s)-c(x,\rho_s))ds}V^\ast(x)\nonumber\\
&\ge& \int_0^{t_1} e^{-\int_0^\tau (q_x(\rho_s)-c(x,\rho_s))ds}\int_{S}V^\ast(y)\tilde{q}(dy|x,\rho_\tau)d\tau +e^{-\int_0^{t_1} (q_x(\rho_s)-c(x,\rho_s))ds}V^\ast(x)
\end{eqnarray}
as follows.

Assume that \begin{eqnarray}\label{ZyExponential018}
\int_0^{t_2} e^{-\int_0^\tau (q_x(\rho_s)-c(x,\rho_s))ds}\int_{S}V^\ast(y)\tilde{q}(dy|x,\rho_\tau)d\tau&<&\infty;\nonumber\\
e^{-\int_0^{t_2} (q_x(\rho_s)-c(x,\rho_s))ds}V^\ast(y)&<&\infty.
\end{eqnarray}
There is no loss of generality in doing so because otherwise (\ref{ZyExponential012}) trivially holds.

Then
\begin{eqnarray}\label{ZyExponential019}
&&\int_0^{t_2} e^{-\int_0^\tau (q_x(\rho_s)-c(x,\rho_s))ds}\int_{S}V^\ast(y)\tilde{q}(dy|x,\rho_\tau)d\tau +e^{-\int_0^{t_2} (q_x(\rho_s)-c(x,\rho_s))ds}V^\ast(x)\nonumber\\
&&-\int_0^{t_1} e^{-\int_0^\tau (q_x(\rho_s)-c(x,\rho_s))ds}\int_{S}V^\ast(y)\tilde{q}(dy|x,\rho_\tau)d\tau \nonumber\\
&&-e^{-\int_0^{t_1} (q_x(\rho_s)-c(x,\rho_s))ds}V^\ast(x)\nonumber\\
&=&\int_{t_1}^{t_2}e^{-\int_0^\tau (q_x(\rho_s)-c(x,\rho_s))ds}\int_{S}V^\ast(y)\tilde{q}(dy|x,\rho_\tau)d\tau\nonumber\\
&&+e^{-\int_0^{t_1} (q_x(\rho_s)-c(x,\rho_s))ds}\left(e^{-\int_{t_1}^{t_2} (q_x(\rho_s)-c(x,\rho_s))ds}-1\right)V^\ast(x)\nonumber\\
&=&\int_0^{t_2-t_1} e^{-\int_0^{t_1+\tau}(q_x(\rho_s)-c(x,\rho_s))ds}\int_S V^\ast(y)\tilde{q}(dy|x,\rho_{t_1+\tau})d\tau\nonumber\\
&&+e^{-\int_0^{t_1} (q_x(\rho_s)-c(x,\rho_s))ds}\left(e^{-\int_{t_1}^{t_2} (q_x(\rho_s)-c(x,\rho_s))ds}-1\right)V^\ast(x)\nonumber\\
&=&\int_0^{t_2-t_1}e^{-\int_0^{t_1} (q_x(\rho_s)-c(x,\rho_s))ds}e^{-\int_{t_1}^{t_1+\tau}(q_x(\rho_s)-c(x,\rho_s))ds}\int_SV^\ast(y)\tilde{q}(dy|x,\rho_{t_1+\tau})d\tau\nonumber\\
&&+e^{-\int_0^{t_1} (q_x(\rho_s)-c(x,\rho_s))ds}\left(e^{-\int_{t_1}^{t_2} (q_x(\rho_s)-c(x,\rho_s))ds}-1\right)V^\ast(x)\nonumber\\
&=&e^{-\int_0^{t_1} (q_x(\rho_s)-c(x,\rho_s))ds}\left\{\int_0^{t_2-t_1} e^{-\int_0^\tau(q_x(\rho_{s+t_1})-c(x,\rho_{s+t_1}))ds}\right.\nonumber\\
&&\left.\times\int_S V^\ast(y)\tilde{q}(dy|x,\rho_{t_1+\tau})d\tau+e^{-\int_{t_1}^{t_2} (q_x(\rho_s)-c(x,\rho_s))ds}V^\ast(x)-V^\ast(x)\right\}.
\end{eqnarray}

Let $\delta>0$ be arbitrarily fixed. By (\ref{ZyExponential015}), there exists some $\hat{\rho}\in{\cal R}$ such that
\begin{eqnarray}\label{ZyExponential017}
V^\ast(x)+\delta&\ge& \int_0^\infty \int_S V^\ast(y)\tilde{q}(dy|x,\hat{\rho}_\tau)e^{-\int_0^\tau (q_x(\hat{\rho}_s)-c(x,\hat{\rho}_s))ds}d\tau+e^{-\int_0^\infty q_x(\hat{\rho}_s)ds}e^{\int_0^\infty c(x,\hat{\rho}_s)ds}.
\end{eqnarray}

Define $\tilde{\rho}\in {\cal R}$ by
\begin{eqnarray}\label{ZyExponential016}
\tilde{\rho}_s=\left\{\begin{array}{ll}
\rho_{t_1+s}, & \mbox{ if }s\le t_2-t_1;  \\
\hat{\rho}_{s-(t_2-t_1)} & \mbox{ if } s>t_2-t_1. \end{array}\right.
\end{eqnarray}
Then
\begin{eqnarray*}
V^\ast(x)&\le& \int_0^\infty e^{-\int_0^\tau (q_x(\tilde{\rho}_s)-c(x,\tilde{\rho}_s))ds} \left(\int_S V^\ast(y)\tilde{q}(dy|x,\tilde{\rho}_\tau)\right)d\tau +e^{-\int_0^\infty q_x(\tilde{\rho}_s)ds}e^{\int_0^\infty c(x,\tilde{\rho}_s)ds}\\
&=&\int_0^{t_2-t_1}e^{-\int_0^\tau (q_x(\rho_{s+t_1})-c(x,\rho_{s+t_1}))ds}\int_S V^\ast(y)\tilde{q}(dy|x,\rho_{t_1+\tau})d\tau\\
&&+\int_{t_2-t_1}^\infty e^{-\int_0^\tau (q_x(\tilde{\rho}_s)-c(x,\tilde{\rho}_s))ds} \left(\int_S V^\ast(y)\tilde{q}(dy|x,\tilde{\rho}_\tau)\right)d\tau +e^{-\int_0^\infty q_x(\tilde{\rho}_s)ds}e^{\int_0^\infty c(x,\tilde{\rho}_s)ds}\\
&=&\int_0^{t_2-t_1}e^{-\int_0^\tau (q_x(\rho_{s+t_1})-c(x,\rho_{s+t_1}))ds}\int_S V^\ast(y)\tilde{q}(dy|x,\rho_{t_1+\tau})d\tau\\
&&+\int_0^\infty  e^{-\int_0^{\tau+(t_2-t_1)}(q_x(\tilde{\rho}_{s})-c(x,\tilde{\rho}_s))ds}\int_S V^\ast(y)\tilde{q}(dy|x,\tilde{\rho}_{\tau+t_2-t_1})d\tau+e^{-\int_0^\infty q_x(\tilde{\rho}_s)ds}e^{\int_0^\infty c(x,\tilde{\rho}_s)ds}\\
&=&\int_0^{t_2-t_1}e^{-\int_0^\tau (q_x(\rho_{s+t_1})-c(x,\rho_{s+t_1}))ds}\int_S V^\ast(y)\tilde{q}(dy|x,\rho_{t_1+\tau})d\tau\\
&&+\int_0^\infty  e^{-\int_0^{t_2-t_1}(q_x(\tilde{\rho}_{s})-c(x,\tilde{\rho}_s))ds}   e^{-\int_{t_2-t_1}^{\tau+(t_2-t_1)}(q_x(\tilde{\rho}_{s})-c(x,\tilde{\rho}_s))ds}  \int_S V^\ast(y)\tilde{q}(dy|x,\tilde{\rho}_{\tau+t_2-t_1})d\tau\\
&&+e^{-\int_0^\infty q_x(\tilde{\rho}_s)ds}e^{\int_0^\infty c(x,\tilde{\rho}_s)ds},
\end{eqnarray*}
where the first inequality is by (\ref{ZyExponential015}). Substituting (\ref{ZyExponential016}) in the last expression, we see
\begin{eqnarray*}
V^\ast(x)&\le& \int_0^{t_2-t_1}e^{-\int_0^\tau (q_x(\rho_{s+t_1})-c(x,\rho_{s+t_1}))ds}\int_S V^\ast(y)\tilde{q}(dy|x,\rho_{t_1+\tau})d\tau\\
&&+e^{-\int_0^{t_2-t_1}    (q_x(\rho_{s+t_1})-c(x,\rho_{s+t_1}))ds}  \int_0^\infty    e^{-\int_{t_{2}-t_{1}}^{\tau+(t_{2}-t_{1})} (q_x(\hat{\rho}_{s-(t_2-t_1)})-c(x,\hat{\rho}_{s-(t_2-t_1)}))ds}\\
&&\times\int_S V^\ast(y)\tilde{q}(dy|x,\hat{\rho}_{\tau})d\tau+e^{-\int_0^{t_2-t_1}(q_x(\rho_{t_1+s})-c(x,\rho_{t_1+s}))ds}e^{-\int_{t_2-t_1}^\infty q_x(\hat{\rho}_{s-(t_2-t_1)})ds}\\
&&\times e^{\int_{t_2-t_1}^\infty c(x,\hat{\rho}_{s-(t_2-t_1)})ds}\\
&=&\int_0^{t_2-t_1}e^{-\int_0^\tau (q_x(\rho_{s+t_1})-c(x,\rho_{s+t_1}))ds}\int_S V^\ast(y)\tilde{q}(dy|x,\rho_{t_1+\tau})d\tau\\
&&+e^{-\int_0^{t_2-t_1}(q_x(\rho_{t_1+s})-c(x,\rho_{t_1+s}))ds}\\
&&\times \left\{\int_0^\infty e^{-\int_0^\tau(q_x(\hat{\rho}_s)-c(x,\hat{\rho}_s))ds}\int_S V^\ast(y)\tilde{q}(dy|x,\hat{\rho}_\tau)d\tau+e^{-\int_0^\infty q_x(\hat{\rho}_s)ds}e^{\int_0^\infty c(x,\hat{\rho}_s)ds}\right\}\\
&\le&\int_0^{t_2-t_1}e^{-\int_0^\tau (q_x(\rho_{s+t_1})-c(x,\rho_{s+t_1}))ds}\int_S V^\ast(y)\tilde{q}(dy|x,\rho_{t_1+\tau})d\tau\\
&&+e^{-\int_0^{t_2-t_1}(q_x(\rho_{t_1+s})-c(x,\rho_{t_1+s}))ds}V^\ast(x)+e^{-\int_0^{t_2-t_1}(q_x(\rho_{t_1+s})-c(x,\rho_{t_1+s}))ds}\delta,
\end{eqnarray*}
where the last inequality is by (\ref{ZyExponential017}). Since $\delta>0$ is arbitrarily fixed, and keeping in mind
(\ref{ZyExponential018}), this amounts to
\begin{eqnarray*}
V^\ast(x)&\le& \int_0^{t_2-t_1}e^{-\int_0^\tau (q_x(\rho_{s+t_1})-c(x,\rho_{s+t_1}))ds}\int_S V^\ast(y)\tilde{q}(dy|x,\rho_{t_1+\tau})d\tau\\
&&+e^{-\int_0^{t_2-t_1}(q_x(\rho_{t_1+s})-c(x,\rho_{t_1+s}))ds}V^\ast(x).
\end{eqnarray*}
This, (\ref{ZyExponential018}) and (\ref{ZyExponential019}) imply (\ref{ZyExponential012}). $\hfill\Box$
\bigskip

\begin{lemma}\label{ZyExponentialLemma05L}
For each $t\ge 0$ and $x\in S$,
\begin{eqnarray*}
\inf_{\rho\in{\cal  R}}\left\{\int_0^t e^{-\int_0^s(q_x(\rho_v)-c(x,\rho_v))dv}\int_S V^\ast(y)\tilde{q}(dy|x,\rho_s)ds+e^{-\int_0^t (q_x(\rho_s)-c(x,\rho_s))ds}V^\ast(x)\right\}=V^\ast(x).
\end{eqnarray*}
\end{lemma}
\par\noindent\textit{Proof.} We only need consider when $t>0$; the case of $t=0$ is trivial. Let $\delta>0$ be arbitrarily fixed. Then by (\ref{ZyExponential015}), there is some $\hat{\rho}\in {\cal R}$ such that
\begin{eqnarray*}
V^\ast(x)+\delta\ge \int_0^\infty e^{-\int_0^\tau (q_x(\hat{\rho}_s)-c(x,\hat{\rho}_s))ds}\int_SV^\ast(y)\tilde{q}(dy|x,\hat{\rho}_\tau)d\tau+e^{-\int_0^\infty q_x(\hat{\rho}_s)ds}e^{-\int_0^\infty c(x,\hat{\rho}_s)ds}.
\end{eqnarray*}
Define $\tilde{\rho}\in {\cal R}$ by
\begin{eqnarray*}
\tilde{\rho}_s=\hat{\rho}_{t+s},~\forall~s>0.
\end{eqnarray*}
Direct calculations similar to those in the proof of Lemma \ref{ZyExponentialLemma02L} show
\begin{eqnarray*}
V^\ast(x)+\delta&\ge& \int_0^t e^{-\int_0^\tau (q_x(\hat{\rho}_s)-c(x,\hat{\rho}_s))ds}\int_S V^\ast(y)\tilde{q}(dy|x,\hat{\rho}_\tau)d\tau+e^{-\int_0^t (q_x(\hat{\rho}_s)-c(x,\hat{\rho}_s))ds}\\
&&\times\left\{\int_0^\infty e^{-\int_0^\tau (q_x(\tilde{\rho}_s)-c(x,\tilde{\rho}_s))ds}\int_S V^\ast(y)\tilde{q}(dy|x,\tilde{\rho}_\tau)d\tau+e^{-\int_0^\infty q_x(\tilde{\rho}_s)ds}e^{-\int_0^\infty c(x,\tilde{\rho}_s)ds}\right\}\\
&\ge&\int_0^t e^{-\int_0^\tau (q_x(\hat{\rho}_s)-c(x,\hat{\rho}_s))ds}\int_S V^\ast(y)\tilde{q}(dy|x,\hat{\rho}_\tau)d\tau+e^{-\int_0^t (q_x(\hat{\rho}_s)-c(x,\hat{\rho}_s))ds}V^\ast(x)
\end{eqnarray*}
where the last inequality is by (\ref{ZyExponential015}). Since $\delta>0$ is arbitrarily fixed, the above implies
\begin{eqnarray*}
V^\ast(x)\ge \inf_{\rho\in {\cal R}}\left\{\int_0^t e^{-\int_0^\tau (q_x({\rho}_s)-c(x,{\rho}_s))ds}\int_S V^\ast(y){q}(dy|x,{\rho}_\tau)d\tau+e^{-\int_0^t (q_x({\rho}_s)-c(x,{\rho}_s))ds}V^\ast(x)\right\}.
\end{eqnarray*}
On the other hand, Lemma \ref{ZyExponentialLemma02L} implies
\begin{eqnarray*}
V^\ast(x)&\le& \inf_{\rho\in {\cal R}}\left\{\int_0^t e^{-\int_0^\tau (q_x(\rho_s)-c(x,\rho_s))ds}\int_{S}V^\ast(y)\tilde{q}(dy|x,\rho_\tau)d\tau +e^{-\int_0^t (q_x(\rho_s)-c(x,\rho_s))ds}V^\ast(x)\right\}.
\end{eqnarray*}
The statement follows from this and the previous inequality. $\hfill\Box$
\bigskip

Under extra conditions, including that the function $q_x(a)-c(x,a)$ is bounded, and $0<\delta<q_x(a)-c(x,a)$ for some constant $\delta>0$, as in \cite{Davis:1993}, the minimization problem on the right hand side of (\ref{ZyExponential015}) can be reduced to a problem of Mayer form, and then
Lemmas \ref{ZyExponentialLemma02L} and \ref{ZyExponentialLemma05L} follow from Lemma (45.12) of \cite{Davis:1993}.

For the next two lemmas, it is convenient to introduce the following notation.
For each $x\in S$ and $T>0,$ let ${\cal R}_{V^\ast,x,T}$ be the set of $\rho\in {\cal R}$ such that
\begin{eqnarray*}
\int_0^te^{-\int_0^s (q_x(\rho_v)-c(x,\rho_v))dv}\int_S V^\ast(y)\tilde{q}(dy|x,\rho_s)ds+e^{-\int_0^t (q_x(\rho_s)-c(x,\rho_s))ds}V^\ast(x)<\infty,~\forall~t\in(0,T).
\end{eqnarray*}

Since $V^\ast(x)\ge 1,$ the above inequality is equivalent to
\begin{eqnarray}\label{ZyExponential020}
&&\int_0^te^{-\int_0^s (q_x(\rho_v)-c(x,\rho_v))dv}\int_S V^\ast(y)\tilde{q}(dy|x,\rho_s)ds<\infty,~e^{-\int_0^t (q_x(\rho_s)-c(x,\rho_s))ds}<\infty,~\forall~t\in(0,T);\nonumber\\
&&V^\ast(x)<\infty,
\end{eqnarray}
for each $\rho\in {\cal R}_{V^\ast,x,T}.$ Note that if $x\in S$ is such that $V^\ast(x)<\infty,$ then by Lemmas \ref{ZyExponentialLemma02L} and \ref{ZyExponentialLemma05L},  ${\cal R}_{V^\ast,x,T}\ne \emptyset.$

\begin{lemma}\label{ZyExponential023}
Let $x\in S$ and $T>0$ be fixed. For each $\rho\in {\cal R}_{V^\ast,x,T}$, it holds that
\begin{eqnarray*}
\int_S V^\ast(y)\tilde{q}(dy|x,\rho_s)\ge V^\ast(x)(q_x(\rho_s)-c(x,\rho_s))
\end{eqnarray*}
almost everywhere with respect to $s\in(0,T).$
\end{lemma}
\par\noindent\textit{Proof.} Since $\rho\in  {\cal R}_{V^\ast,x,T}$, one can apply the fundamental theorem of calculus and differentiate
\begin{eqnarray*}
\int_0^t e^{-\int_0^\tau (q_x(\rho_s)-c(x,\rho_s))ds}\int_{S}V^\ast(y)\tilde{q}(dy|x,\rho_\tau)d\tau +e^{-\int_0^t (q_x(\rho_s)-c(x,\rho_s))ds}V^\ast(x)
\end{eqnarray*}
with respect to $t\in(0,T),$ and deduce
\begin{eqnarray*}
e^{-\int_0^t (q_x(\rho_s)-c(x,\rho_s))ds}\int_S V^\ast(y)\tilde{q}(dy|x,\rho_t)-e^{-\int_0^t(q_x(\rho_s)-c(x,\rho_s))ds}(q_x(\rho_t)-c(x,\rho_t))V^\ast(x)\ge 0
\end{eqnarray*}
for almost all $t\in(0,T)$ with respect to the Lebesgue measure, where the last inequality is by Lemma \ref{ZyExponentialLemma02L}.
The statement of the lemma immediately follows. (Recall that (\ref{ZyExponential020}) holds for each $\rho\in {\cal R}_{V^\ast,x,T}$.) $\hfill\Box$
\bigskip

\begin{lemma}\label{ZyExponentialLemma06L}
For each $x\in S$, where $V^\ast(x)<\infty$, (\ref{ZyExponential036}) is satisfied.
\end{lemma}
\par\noindent\textit{Proof.} Let $T>0$ be arbitrarily fixed, and so is $x\in S$, where $V^\ast(x)<\infty$. Then ${\cal R}_{V^\ast,x,T}\ne \emptyset$ as explained earlier. Let some $\rho\in {\cal R}_{V^\ast,x,T}$ be arbitrarily fixed. One can legitimately write
\begin{eqnarray*}
e^{-\int_0^t(q_x(\rho_s)-c(x,\rho_s))ds}V^\ast(x)-V^\ast(x)&=&-\int_{0}^t (q_x(\rho_\tau)-c(x,\rho_\tau))e^{-\int_0^\tau(q_x(\rho_s)-c(x,\rho_s))ds}d\tau V^\ast(x),\\
&&~\forall~t\in(0,T).
\end{eqnarray*}
Now,
\begin{eqnarray*}
&&\int_0^t e^{-\int_0^s(q_x(\rho_v)-c(x,\rho_v))dv}\int_S V^\ast(y)\tilde{q}(dy|x,\rho_s)ds+e^{-\int_0^t (q_x(\rho_s)-c(x,\rho_s))ds}V^\ast(x)-V^\ast(x)\\
&=&\int_0^t e^{-\int_0^\tau(q_x(\rho_v)-c(x,\rho_v))dv}\int_S V^\ast(y)\tilde{q}(dy|x,\rho_\tau)d\tau\\
&&-\int_{0}^t (q_x(\rho_\tau)-c(x,\rho_\tau))e^{-\int_0^\tau(q_x(\rho_s)-c(x,\rho_s))ds}d\tau V^\ast(x)\\
&=&\int_0^t e^{-\int_0^\tau (q_x(\rho_s)-c(x,\rho_s))ds}\left\{\int_S V^\ast(y)\tilde{q}(dy|x,\rho_\tau)-(q_x(\rho_\tau)-c(x,\rho_\tau))V^\ast(x)\right\}d\tau\\
&=&\int_0^t e^{-\int_0^\tau (q_x(\rho_s)-c(x,\rho_s))ds}\int_A\rho_\tau(da)\left\{\int_S V^\ast(y)\tilde{q}(dy|x,a)-(q_x(a)-c(x,a))V^\ast(x)\right\}d\tau
\end{eqnarray*}
for each $t\in(0,T).$

By Lemma \ref{ZyExponentialLemma05L}, we deduce from the above that
\begin{eqnarray}\label{ZyExponential021}
&&0\nonumber\\
&=&\inf_{\rho\in {\cal R}_{V^\ast,x,T}}\left\{\int_0^t e^{-\int_0^\tau (q_x(\rho_s)-c(x,\rho_s))ds}\int_A\rho_\tau(da)\left\{\int_S V^\ast(y)\tilde{q}(dy|x,a)-(q_x(a)-c(x,a))V^\ast(x)\right\}d\tau\right\}\nonumber\\
&\ge&\inf_{\rho\in {\cal R}_{V^\ast,x,T}}\left\{\int_0^t e^{-\int_0^\tau (q_x(\rho_s)-c(x,\rho_s))ds}\inf_{a\in A}\left\{\int_S V^\ast(y)\tilde{q}(dy|x,a)-(q_x(a)-c(x,a))V^\ast(x)\right\}d\tau\right\}\nonumber\\
&\ge&\inf_{\rho\in {\cal R}_{V^\ast,x,T}}\left\{\int_0^t e^{-\tau \overline{q}_x}\inf_{a\in A}\left\{\int_S V^\ast(y)\tilde{q}(dy|x,a)-(q_x(a)-c(x,a))V^\ast(x)\right\}d\tau\right\}
\end{eqnarray}
where the first equality is also because of $V^\ast(x)<\infty.$
Let \begin{eqnarray*}
B(x)=\left\{a\in A:~ \int_S V^\ast(y)\tilde{q}(dy|x,a)<\infty\right\}.
\end{eqnarray*}
Then
\begin{eqnarray}\label{ZyExponential022}
&&\inf_{a\in A}\left\{\int_S V^\ast(y)\tilde{q}(dy|x,a)-(q_x(a)-c(x,a))V^\ast(x)\right\}\nonumber\\
&=&\inf_{a\in B(x)}\left\{\int_S V^\ast(y)\tilde{q}(dy|x,a)-(q_x(a)-c(x,a))V^\ast(x)\right\}
\end{eqnarray}

Each $b\in B(x)$ is identified by an element $\rho^b\in {\cal R}$ such that $\rho^b_s(da)=\delta_{\{b\}}(da)$ for all $s>0.$ Furthermore, compared with (\ref{ZyExponential020}) and keeping in mind $V^\ast(x)<\infty,$ we see that $\rho^b\in {\cal R}_{V^\ast,x,T}.$ Thus, $\{\rho^b:~b\in B(x)\}\subseteq  {\cal R}_{V^\ast,x,T}.$ Hence, for each $a\in B(x),$ one can apply Lemma \ref{ZyExponential023}, and after that, see
\begin{eqnarray*}
\inf_{a\in B(x)}\left\{\int_S V^\ast(y)\tilde{q}(dy|x,a)-(q_x(a)-c(x,a))V^\ast(x)\right\}\ge 0.
\end{eqnarray*}
Consequently, we see from (\ref{ZyExponential021}), (\ref{ZyExponential022}) and the above inequality that
\begin{eqnarray*}
0\ge \inf_{\rho\in {\cal R}_{V^\ast,x,T}}\left\{\int_0^t e^{-\tau \overline{q}_x}\inf_{a\in B(x)}\left\{\int_S V^\ast(y)\tilde{q}(dy|x,a)-(q_x(a)-c(x,a))V^\ast(x)\right\}d\tau\right\}\ge 0,
\end{eqnarray*}
and thus
\begin{eqnarray*}
0=\inf_{a\in A}\left\{c(x,a)V^\ast(x)+\int_S q(dy|x,a)V^\ast(y)\right\},
\end{eqnarray*}
as required. (Recall that $\overline{q}_x<\infty$.) $\hfill\Box$
\bigskip

Now we are ready to present the proof of Theorem \ref{ZyExponentialTheorem06} as follows.
\bigskip

\par\noindent\textit{Proof of Theorem \ref{ZyExponentialTheorem06}.} Part (a) of this statement has been proved in Lemma \ref{ZyExponentialLemma06L}. We prove part (b) of the statement as follows.
For each $x\in S,$ we can view $\varphi(x)$ as an element of ${\cal R}$ by identifying it with $\rho^x\in {\cal R}$ such that
\begin{eqnarray*}
\rho^x_t(da):=\delta_{\{\varphi(x)\}}(da),~\forall~ t>0.
\end{eqnarray*}
Then, $x\in S\rightarrow \rho^x$ clearly defines a specific deterministic stationary strategy for the DTMDP model $\{\textbf{X},\textbf{A},p,l\}$ defined in Section \ref{ZyExponentialSec02}.

Let $x\in S$ be arbitrarily fixed, where $V^\ast(x)<\infty.$ Then by (\ref{ZyExponential026}),
\begin{eqnarray}\label{ZyExponential028}
(q_x(\varphi(x))-c(x,\varphi(x)))V^\ast(x)=\int_S \tilde{q}(dy|x,\varphi(x))V^\ast(y).
\end{eqnarray}
Note that the right hand side is nonnegative, and is zero if and only if $q_x(\varphi(x))=0$, because $V^\ast(x)\ge 1$ for each $x\in S.$
In other words, there are only two possibilities;
\begin{eqnarray}\label{ZyExponential027}
q_x(\varphi(x))> c(x,\varphi(x))\ge 0,
\end{eqnarray}
or
\begin{eqnarray}\label{ZyExponential029}
q_x(\varphi(x))=c(x,\varphi(x))=0.
\end{eqnarray}

In case of (\ref{ZyExponential027}), we see
\begin{eqnarray*}
&&V^\ast(x)=\inf_{\rho\in {\cal R}}\left\{\int_0^\infty e^{-\int_0^\tau (q_x(\rho_s))-c(x,\rho_s))ds} \left(\int_S V^\ast(y)\tilde{q}(dy|x,\rho_\tau)\right)d\tau +e^{-\int_0^\infty q_x(\rho_s)ds}e^{\int_0^\infty c(x,\rho_s)ds} \right\}\\
&\le& \int_0^\infty e^{-\int_0^\tau (q_x(\rho^x_s))-c(x,\rho^x_s))ds} \left(\int_S V^\ast(y)\tilde{q}(dy|x,\rho^x_\tau)\right)d\tau +e^{-\int_0^\infty q_x(\rho^x_s)ds}e^{\int_0^\infty c(x,\rho^x_s)ds}\\
&=&\int_0^\infty  e^{-(q_x(\varphi(x))-c(x,\varphi(x)))\tau}(q_x(\varphi(x))-c(x,\varphi(x)))d\tau V^\ast(x)+e^{-\int_0^\infty q_x(\varphi(x))ds}e^{\int_0^\infty c(x,\varphi(x))ds}\\
&=& V^\ast(x),
\end{eqnarray*}
where the first equality is by (\ref{ZyExponential015}), and the second equality is by (\ref{ZyExponential028}) and the definition of $\rho^x$, and the last equality is by (\ref{ZyExponential027}); recall (\ref{ZyExponential56}). Thus,
\begin{eqnarray*}
&&\inf_{\rho\in {\cal R}}\left\{\int_0^\infty e^{-\int_0^\tau (q_x(\rho_s))-c(x,\rho_s))ds} \left(\int_S V^\ast(y)\tilde{q}(dy|x,\rho_\tau)\right)d\tau +e^{-\int_0^\infty q_x(\rho_s)ds}e^{\int_0^\infty c(x,\rho_s)ds} \right\}\\
&=& \int_0^\infty e^{-\int_0^\tau (q_x(\rho^x_s))-c(x,\rho^x_s))ds} \left(\int_S V^\ast(y)\tilde{q}(dy|x,\rho^x_\tau)\right)d\tau +e^{-\int_0^\infty q_x(\rho^x_s)ds}e^{\int_0^\infty c(x,\rho^x_s)ds}
\end{eqnarray*}
under (\ref{ZyExponential027}). Similar calculation show that in case of (\ref{ZyExponential029}), and in case of $V^\ast(x)=\infty,$ the above equalities hold as well. It remains to apply Proposition \ref{ZyExponentialProposition01}(b); recall the discussions in Section \ref{ZyExponentialSec02} about the reduction of the CTMDP model $\{S,A,q,c\}$ to the DTMDP model $\{\textbf{X},\textbf{A},p,l\}$ therein. $\hfill\Box$

\subsection{Existence of a deterministic stationary optimal policy}
The objective of this subsection is to show the existence of a deterministic stationary optimal policy for the CTMDP problem (\ref{ZyExponChap1}), under some compactness-continuity conditions.

From now on, the following assumption is always in place.
\begin{assumption}
For each $x\in S,$ \begin{eqnarray*}
\sup_{a\in A}\{c(x,a)\}=:\overline{c}(x)<\infty.
\end{eqnarray*}
\end{assumption}
We mention that the function $\overline{c}$ is upper semianalytic on $S$, and may be not Borel measurable; the similar remark holds for the function $\overline{q}$, see \cite{Bertsekas:1978}. However, we have the following handy fact\footnote{I was told the fact in Lemma \ref{ZyExponentialLemma07L} by Professor Eugene A. Feinberg.}.

\begin{lemma}\label{ZyExponentialLemma07L}
There exists a $[1,\infty)$-valued Borel measurable function $w$ on $S$ such that
\begin{eqnarray}\label{ZyExponential50}
w(x)\ge 1+\overline{c}(x)+\overline{q}_x,~\forall~x\in S.
\end{eqnarray}
\end{lemma}
\par\noindent\textit{Proof.} This follows from the reasoning of the proof of Lemma 1(a) in \cite{Feinberg:2016} based on the Novikov seperation theorem. $\hfill\Box$\bigskip

The role of the function $w$ can be also well appreciated in the next subsection.

Next we present two sets of compactness-continuity conditions, under either of which, the main optimality results presented henceforth survive.
\begin{condition}\label{ZyExponentialCondition01}
\begin{itemize}
\item[(a)] The function $w$ from Lemma \ref{ZyExponentialLemma07L} is continuous on $S.$
\item[(b)] For each bounded continuous function $f$ on $S$, $\int_S f(y)\tilde{q}(dy|x,a)$ is continuous in $(x,a)\in S\times A.$
\item[(c)] The function $c(x,a)$ is lower semicontinuous in $(x,a)\in S\times A.$
\item[(d)] The action space $A$ is a compact Borel space.
\end{itemize}
\end{condition}
Part (a) of Condition \ref{ZyExponentialCondition01} is not restrictive, see p.48 of \cite{Srivastava:1998}.

\begin{condition}\label{ZyExponentialCondition02}
\begin{itemize}
\item[(a)] For each bounded measurable function $f$ on $S$ and each $x\in S$, $\int_S f(y)\tilde{q}(dy|x,a)$ is continuous in $a\in A.$
\item[(b)] For each $x\in S,$ the function $c(x,a)$ is lower semicontinuous in $a\in A.$
\item[(c)] The action space $A$ is a compact Borel space.
\end{itemize}
\end{condition}

Condition \ref{ZyExponentialCondition01} is called the compactness-weak continuity condition, and Condition \ref{ZyExponentialCondition02} is called the compactness-strong continuity condition. Often, the weak continuity condition is easier for verifications, and it is noted that in some practical applications, the weak continuity condition is satisfied while the strong continuity condition is not, see e.g., Section 6 of \cite{Jaskiewicz:2009}. Nevertheless, the two conditions do not imply each other.

\begin{theorem}\label{ZyExponentialTheorem07}
 Suppose Condition \ref{ZyExponentialCondition01} (resp., Condition \ref{ZyExponentialCondition02}) is satisfied. Then there exists a deterministic stationary policy $\varphi$ satisfying (\ref{ZyExponential026}) for each $x\in S,$ where $V^\ast(x)<\infty$, and so there exists a deterministic stationary optimal policy $\varphi$ for the CTMDP problem (\ref{ZyExponChap1}), and the value function $V^\ast$ is lower semicontinuous (resp., measurable) on $S$.
\end{theorem}
We postpone the proof of Theorem \ref{ZyExponentialTheorem07} after the next few lemmas and preliminaries.

Let us equip ${\cal R}$ with the Young topology, which is the weakest
topology with respect to which the function
\begin{eqnarray*}
\rho\in {\cal{R}}\rightarrow \int_0^\infty \int_A
f(t,a)\rho_t(da)dt
\end{eqnarray*}
is continuous for each strongly integrable Carath\'{e}odory functions
$f$ on $(0,\infty)\times A$ . Here a real-valued measurable
function $f$ on $(0,\infty)\times A$ is called a strongly
integrable Carath\'{e}odory function if for each fixed
$t\in(0,\infty)$, $f(t,a)$ is continuous in $a\in A,$ and for each
fixed $a\in A,$ $\sup_{a\in A}|f(t,a)|$ is integrable in $t$,
i.e., $\int_0^\infty \sup_{a\in A}|f(t,a)|dt<\infty.$ See more details in \cite{Davis:1993}.

\begin{lemma}\label{ZyExponentialYoungtopologycompact}
Endowed with the Young topology, if the action space $A$ is compact, then ${\cal
R}$ is a compact Borel space.
\end{lemma}
\par\noindent\textit{Proof.} See Remark 8.2.3 of \cite{BauerleRieder:2011}, or Chapter 4 of \cite{Davis:1993}. $\hfill\Box$
\bigskip

\begin{lemma}\label{ZyExponentialLemma08L}
\begin{itemize}
\item[(a)]Suppose Condition \ref{ZyExponentialCondition01} is satisfied. Then the DTMDP model $\{\textbf{X},\textbf{A},p,l\}$ defined in Section \ref{ZyExponentialSec02} satisfies Condition \ref{ZyExponentialDTMDPCon1}.
    \item[(b)] Suppose Condition \ref{ZyExponentialCondition02} is satisfied. Then the DTMDP model $\{\textbf{X},\textbf{A},p,l\}$ defined in Section \ref{ZyExponentialSec02} satisfies Condition \ref{ZyExponentialDTMDPCon2}.
\end{itemize}
\end{lemma}
\par\noindent\textit{Proof.} One can apply the reasoning in the proof of Lemma 3.2 of \cite{Schal:1998}. $\hfill\Box$
\bigskip

\begin{lemma}\label{ZyExponentialLemma10L}
Define the stochastic kernel $\tilde{p}$ on ${\cal B}(S)$ from $(x,a)\in S\times A$ by
\begin{eqnarray*}
\tilde{p}(dy|x,a):=\frac{q(dy|x,a)}{w(x)}+\delta_{\{x\}}(dy),~\forall~(x,a)\in S\times A.
\end{eqnarray*}
Then the following assertions hold.
\begin{itemize}
\item[(a)] An $[1,\infty]$-valued lower semianalytic function $V$ on $S$ satisfies
\begin{eqnarray}\label{ZyExponential037}
0=\inf_{a\in A}\left\{c(x,a)V(x)+\int_S q(dy|x,a)V(y)\right\}
\end{eqnarray}
for each $x\in S$ such that $V(x)<\infty$ if and only if $V$ is an $[1,\infty]$-valued lower semianalytic solution to
\begin{eqnarray}\label{ZyExponential038}
V(x)=\inf_{a\in A}\left\{\frac{w(x)}{w(x)-c(x,a)}\int_{S}\tilde{p}(dy|x,a)V(y)\right\},~\forall~x\in S.
\end{eqnarray}
\item[(b)] Let $V$ be an $[1,\infty]$-valued lower semianalytic function on $S$ satisfying (\ref{ZyExponential037}) for each $x\in S$ such that $V(x)<\infty.$ A deterministic stationary policy $\varphi$ satisfies
\begin{eqnarray}\label{ZyExponential51}
0=\inf_{a\in A}\left\{c(x,a)V(x)+\int_S q(dy|x,a)V(y)\right\}=c(x,\varphi(x))V(x)+\int_S q(dy|x,\varphi(x))V(y)
\end{eqnarray}
for each $x\in S$ such that $V(x)<\infty$ if and only if this deterministic stationary policy $\varphi$ satisfies
    \begin{eqnarray*}
    \inf_{a\in A}\left\{\frac{w(x)}{w(x)-c(x,a)}\int_{S}\tilde{p}(dy|x,a)V(y)\right\}= \frac{w(x)}{w(x)-c(x,\varphi(x))}\int_{S}\tilde{p}(dy|x,\varphi(x))V(y),~\forall~x\in S.
    \end{eqnarray*}
\end{itemize}
\end{lemma}
\par\noindent\textit{Proof.} (a) We first show the ``only if'' part. Let $V$ be a $[1,\infty]$-valued lower semianalytic solution (\ref{ZyExponential037}). Let $x\in S$ be fixed. If $V(x)=\infty,$ then (\ref{ZyExponential038}) is satisfied as the both sides are infinite. Suppose now $V(x)<\infty.$ Then
\begin{eqnarray*}
0&\le& c(x,a)V(x)+w(x)\int_S \frac{q(dy|x,a)}{w(x)}V(y)+w(x)V(x)-w(x)V(x)\\
&=&c(x,a)V(x)+w(x)\int_S \tilde{p}(dy|x,a)V(y)-w(x)V(x),~\forall~a\in A.
\end{eqnarray*}
Following from this and keeping in mind (\ref{ZyExponential50}), simple calculations imply
\begin{eqnarray*}
V(x)\le \inf_{a\in A}\left\{\frac{w(x)}{w(x)-c(x,a)}\int_S \tilde{p}(dy|x,a)V(y)\right\}.
\end{eqnarray*}
To show the equality, let $\epsilon>0$ be arbitrarily fixed. By (\ref{ZyExponential037}), there exists some $\hat{a}\in A$ such that
\begin{eqnarray*}
\epsilon>c(x,\hat{a})V(x)+\int_S q(dy|x,\hat{a})V(y).
\end{eqnarray*}
The above inequality implies
\begin{eqnarray*}
V(x)+ \epsilon \ge V(x)+\frac{\epsilon}{w(x)-c(x,\hat{a})}> \frac{w(x)}{w(x)-c(x,\hat{a})}\int_{S}\tilde{p}(dy|x,\hat{a})V(y),
\end{eqnarray*}
where the first inequality is by (\ref{ZyExponential50}). Thus, (\ref{ZyExponential038}) is satisfied.

The similar reasoning applies to show the ``if'' part; the details are omitted.

(b) This part can be proved as for part (a). $\hfill\Box$
\bigskip

\begin{lemma}\label{ZyExponentialLemma09L}
\begin{itemize}
\item[(a)] Suppose Condition \ref{ZyExponentialCondition01} is satisfied. Then $(x,a)\in S\times A\rightarrow \frac{w(x)}{w(x)-c(x,a)}$ is lower semicontinuous, and for each bounded continuous function $f$ on $S$, $(x,a)\in S\times A\rightarrow \int_{S}f(y)\tilde{p}(dy|x,a)$ is continuous.
\item[(b)] Suppose Condition \ref{ZyExponentialCondition02} is satisfied. Then for each $x\in S,$ $a\in A\rightarrow \frac{w(x)}{w(x)-c(x,a)}$ is lower semicontinuous, and for each bounded measurable function $f$ on $S$, $a\in A\rightarrow \int_{S}f(y)\tilde{p}(dy|x,a)$ is continuous.
\end{itemize}
\end{lemma}
\par\noindent\textit{Proof.} The statement of this lemma is immediate from Condition \ref{ZyExponentialCondition01} and Condition \ref{ZyExponentialCondition02}, respectively, as well as the definition of the stochastic kernel $\tilde{p}.$ $\hfill\Box$
\bigskip

Now we are in position to prove Theorem \ref{ZyExponentialTheorem07}.
\bigskip

\par\noindent\textit{Proof of Theorem \ref{ZyExponentialTheorem07}.} Suppose Condition \ref{ZyExponentialCondition01} is satisfied. By Lemma \ref{ZyExponentialLemma08L}, one can apply Proposition \ref{ZyExponentialProposition02}(a) to the DTMDP model $\{\textbf{X},\textbf{A},p,l\}$ defined in Section \ref{ZyExponentialSec02}. This and Theorem \ref{ZyExponentialTheorem05} imply that the value function $V^\ast$ of the CTMDP problem (\ref{ZyExponChap1}) is $[1,\infty]$-valued and lower semicontinuous on $S$. By Theorem \ref{ZyExponentialTheorem06}(a), Lemmas \ref{ZyExponentialLemma10L} and \ref{ZyExponentialLemma09L}, Proposition 7.31 of \cite{Bertsekas:1978} and a well known measurable selection theorem, see e.g., Proposition 7.33 of \cite{Bertsekas:1978}, we see that there is a deterministic stationary policy $\varphi$ for the CTMDP model such that (\ref{ZyExponential026}) is satisfied for each $x\in S$, where $V^\ast(x)<\infty.$ By Theorem \ref{ZyExponentialTheorem06}(b), this deterministic stationary policy $\varphi$ is optimal for the CTMDP problem (\ref{ZyExponChap1}).

The case when Condition \ref{ZyExponentialCondition02} is satisfied can be proved in the same way, by applying Proposition \ref{ZyExponentialProposition02}(b) and the corresponding measurable selection theorem, see e.g., Proposition D.5 of \cite{Hernandez-Lerma:1996}.
$\hfill\Box$\bigskip

\subsection{Further reduction to a simpler DTMDP and value iteration}
In this subsection, we reduce the CTMDP model $\{S,A,q,c\}$ to a DTMDP $\{S,A,\tilde{p},\tilde{l}\}$ with the cost function $\tilde{l}$ being defined below. Compared to the DTMDP model $\{\textbf{X},\textbf{A},p,l\}$ defined in Section \ref{ZyExponentialSec02}, the DTMDP model here is simpler, with the same state and action space as the original CTMDP model.

\begin{theorem}\label{ZyExponentialTheorem08}
\begin{itemize}
\item[(a)] Suppose Condition \ref{ZyExponentialCondition01} is satisfied. Then the value function $V^\ast$ for the CTMDP problem (\ref{ZyExponChap1}) is the minimal $[1,\infty]$-valued lower semicontinuous solution to the optimality equation (\ref{ZyExponential036}).

\item[(b)] Suppose Condition \ref{ZyExponentialCondition02} is satisfied. Then the value function $V^\ast$ for the CTMDP problem (\ref{ZyExponChap1}) is the minimal $[1,\infty]$-valued measurable solution to (\ref{ZyExponential036}).
\end{itemize}
\end{theorem}

\par\noindent\textit{Proof.} (a) Let $V$ be a lower semicontinuous $[1,\infty]$-valued function on $S$ such that (\ref{ZyExponential037}) is satisfied wherever $V(x)<\infty.$ By Proposition 7.33 of \cite{Bertsekas:1978}, there exists a deterministic stationary policy $\varphi$ satisfying (\ref{ZyExponential51}) for each $x\in S$ such that $V(x)<\infty.$ Let $x\in S$ be fixed. Assume for now $V(x)<\infty.$ Arguing as in the proof of Theorem \ref{ZyExponentialTheorem06}(b), we see
\begin{eqnarray*}
&&\inf_{\rho\in {\cal R}}\left\{\int_0^\infty e^{-\int_0^\tau (q_x(\rho_s))-c(x,\rho_s))ds} \left(\int_X V(y)\tilde{q}(dy|x,\rho_\tau)\right)d\tau +e^{-\int_0^\infty q_x(\rho_s)ds}e^{\int_0^\infty c(x,\rho_s)ds} \right\}\\
&\le& \int_0^\infty e^{-\tau(q_x(\varphi(x))-c(x,\varphi(x)))}(c(x,\varphi(x))-q_x(\varphi(x)))d\tau V(x)+e^{-\int_0^\infty q_x(\varphi(x))ds}e^{\int_0^\infty c(x,\varphi(x))ds}\\
&\le& V(x),
\end{eqnarray*}
where for the last inequality, recall that one only needs deal with two possibilities, namely, (\ref{ZyExponential027}) and (\ref{ZyExponential029}), because of (\ref{ZyExponential51}). The previous inequality holds trivially if $V(x)=\infty.$ Now by applying
Proposition \ref{ZyExponentialProposition01} to the DTMDP model $\{\textbf{X},\textbf{A},p,l\}$ in Section \ref{ZyExponentialSec02}, c.f., (\ref{ZyExponential015}), as well as Theorem \ref{ZyExponentialTheorem05}, we see $V(x)\ge V^\ast(x)$ for each $x\in S.$ This and Theorem  \ref{ZyExponentialTheorem06}(a) imply part (a) of this statement.

(b) This part can be proved in the same way as for part (a) by using the appropriate measurable selection theorem, c.f. Proposition D.5 of \cite{Hernandez-Lerma:1996}.
$\hfill\Box$
\bigskip

Now we are in position to present the equivalent DTMDP model $\{S,A,\tilde{p},\tilde{l}\}$.

Define for each $(x,a,y)\in S\times A\times S,$
\begin{eqnarray*}
\tilde{l}(x,a,y):=\ln\frac{w(x)}{w(x)-c(x,a)}.
\end{eqnarray*}
Recall by (\ref{ZyExponential50}), $\tilde{l}(x,a,y)>0$ for each $(x,a,y)\in S\times A\times S.$

Consider the DTMDP model $\{S,A,\tilde{p},\tilde{l}\}$ (with the exponential utility). Note that (\ref{ZyExponential038}) is the optimality equation for this DTMDP model. Suppose Condition \ref{ZyExponentialCondition01} or Condition \ref{ZyExponentialCondition02} is satisfied. Then by Theorem \ref{ZyExponentialTheorem08}, Proposition \ref{ZyExponentialProposition02} applied to the DTMDP model $\{S,A,\tilde{p},\tilde{l}\}$, and Lemma \ref{ZyExponentialLemma10L}(b), we see the DTMDP model $\{S,A,\tilde{p},\tilde{l}\}$ and the CTMDP model $\{S,A,q,c\}$ are equivalent; the value functions are the same, and a deterministic stationary optimal policy for the CTMDP model gives a deterministic stationary optimal strategy for the DTMDP model $\{S,A,\tilde{p},\tilde{l}\}$, and vice versa.

As a consequence, we can write down the value iteration algorithm for the CTMDP problem (\ref{ZyExponChap1}).
\begin{corollary}
Suppose Condition \ref{ZyExponentialCondition01} or Condition \ref{ZyExponentialCondition02} is satisfied. Define $V^{(0)}(x):=1$ for each $x\in S,$ and for each $n=1,2,\dots,$
\begin{eqnarray*}
V^{(n)}(x)=\inf_{a\in A}\left\{\int_S \tilde{p}(dy|x,a)e^{\tilde{l}(x,a,y)}V^{(n-1)}(y)\right\},~\forall~x\in S.
\end{eqnarray*}
Then for each $x\in  S,$ $V^{(n)}(x)$ increases to $V^\ast(x)$ as $n\uparrow \infty$, where $V^\ast$ is the value function of the CTMDP problem (\ref{ZyExponChap1}).
\end{corollary}
\par\noindent\textit{Proof.} The statement follows from the discussions above this theorem, and Proposition \ref{ZyExponentialProposition02}. Recall Lemma \ref{ZyExponentialLemma09L}. $\hfill\Box$

\section{Conclusion}\label{ZyExponentialSect04}
To sum up, for the CTMDP problem, where the certainty equivalent with respect to the exponential utility of the total undiscounted cost is to be minimized, we established the optimality equation. Under the compactness-continuity condition, we showed the existence of a deterministic stationary optimal policy. By investigating the optimality equation, we reduced the CTMDP problem to an equivalent DTMDP problem, which is with the same state and action space as the original CTMDP. In particular, the value iteration algorithm for the CTMDP problem follows from this reduction. Note that, we did not need impose any condition on the growth of the transition rate, and the cost rate is unbounded, and the controlled process in the Borel state space could be explosive.

As for applications, we believe that our results will be useful for optimal control of queueing systems. In fact, \cite{Coraluppi:1997} considered the risk-sensitive control of a queueing system as a CTMDP with the total discounted cost in a finite state and action space. Restricted to deterministic stationary policies, although the authors of \cite{Coraluppi:1997} applied the uniformization technique to reduce the CTMDP to a DTMDP, they did not manage to investigate the continuous-time problem because the resulting DTMDP was nonstandard. The reduction method here is different from the uniformization technique, and its application to the discounted problem will be more delicate.
\bigskip

\par\noindent\textbf{Acknowledgement.} This work was carried out with a financial grant from the Research Fund for Coal and Steel of the European Commission, within the INDUSE-2-SAFETY project (Grant No. RFSR-CT-2014-00025). I would like to thank Alexey B. Piunovskiy for the helpful discussions. Finally, I would like to thank the referees for the careful reading and useful remarks.

\end{document}